\documentclass[a4paper,conference]{IEEEtran}
\IEEEoverridecommandlockouts
\usepackage{cite}
\usepackage{array}
\usepackage{graphicx}
\usepackage{amsmath,amsfonts,amssymb,amsthm}
\usepackage{mathtools}
\usepackage{xcolor}

\usepackage{float}
\usepackage[caption=false,font=footnotesize,labelfont=sf,textfont=sf]{subfig}
\usepackage{tikz}
\usepackage{enumitem}
\usepackage{textcomp}
\usepackage{pgfplots}
\pgfplotsset{compat=1.3}
\pgfplotsset{
    discard if/.style 2 args={
        x filter/.code={
            \edef\tempa{\thisrow{#1}}
            \edef\tempb{#2}
            \ifx\tempa\tempb
                
            \fi
        }
    },
    discard if not/.style 2 args={
        x filter/.code={
            \edef\tempa{\thisrow{#1}}
            \edef\tempb{#2}
            \ifx\tempa\tempb
            \else
                
            \fi
        }
    }
}

\newcommand{\Field}[1]{\mathbb{#1}}
\newcommand{\Measure}[1]{\mathbb{#1}}
\newcommand{\Space}[1]{\mathcal{#1}}
\newcommand{\Set}[1]{\mathcal{#1}}

\newcommand{\Matrix}[1]{\mathbf{#1}}
\newcommand{\Vector}[1]{\pmb{#1}}
\newcommand{\Operator}[1]{\mathsf{#1}}
\newcommand{\Functional}[1]{\mathtt{#1}}

\newcommand{\abs}[1]{|#1|}
\newcommand{\norm}[1]{\lVert#1\rVert}

\newcommand{\lev}{{\mathrm{lev}}}
\newcommand{\innerprod}[2]{\left\langle#1, #2\right\rangle}

\newcommand{\Exp}[1]{\mathbb{E}[ #1 ]}
\newcommand{\Expp}[2]{\mathbb{E}_{#1}\left[ #2 \right]}

\newcommand{\probMeas}{\Measure{P}}

\newcommand{\est}[1]{\hat{#1}}

\newcommand{\diagm}{\text{Diag}}

\newcommand{\mymark}[2]{\stackrel{\mathclap{\normalfont\mbox{\tiny\sffamily #1}}}{#2}}
\newcommand{\mylabel}[2]{#2\def\@currentlabel{#2}\label{#1}}

\theoremstyle{plain}
\newtheorem{theorem}{Theorem}
\newtheorem{proposition}{Proposition}
\newtheorem{lemma}{Lemma}

\theoremstyle{definition}
\newtheorem{definition}{Definition}
\newtheorem{assumption}{Assumption}

\theoremstyle{remark}
\newtheorem{remark}{Remark}

\newcommand{\Real}{\Field{R}}
\newcommand{\Complex}{\Field{C}}
\newcommand{\Natural}{\Field{N}}

\newcommand{\Realp}{\Real_{\geq 0}}
\newcommand{\Realpp}{\Real_{>0}}

\newcommand{\sX}{\Set{X}}

\newcommand{\vx}{\Vector{x}}

\newcommand{\sN}{\Set{A}}

\newcommand{\sG}{\Set{G}}
\newcommand{\sGi}{\sG_i}

\newcommand{\sE}{\Set{E}}
\newcommand{\sEi}{\sE_i}
\newcommand{\sS}{\Set{S}}
\newcommand{\sSb}{\sS^{c}}

\newcommand{\sC}{\Set{C}}

\newcommand{\sD}{\Set{D}}

\newcommand{\sQ}{\Set{Q}}
\newcommand{\sR}{\Set{R}}

\newcommand{\sM}{\Set{M}}

\newcommand{\cT}{\Theta_{k,i}}
\newcommand{\cS}{\Theta'_{k,i}}

\newcommand{\vh}{\Vector{h}}
\newcommand{\vhki}{\vh_{k,i}}

\newcommand{\bh}{\vh^\star}

\newcommand{\vOne}{\Vector{1}}
\newcommand{\ve}{\Vector{e}}
\newcommand{\veki}{\ve_{k,i}}
\newcommand{\vn}{\Vector{n}}
\newcommand{\vnki}{\vn_{k,i}}
\newcommand{\vv}{\Vector{v}}

\newcommand{\vni}{\vn_i}

\newcommand{\vy}{\Vector{y}}
\newcommand{\vym}{\vy^{(m)}}

\newcommand{\vu}{\Vector{u}}
\newcommand{\vut}{\Tilde{\vu}}
\newcommand{\vz}{\Vector{z}}
\newcommand{\vzki}{\vz_{k,i}}
\newcommand{\vzero}{\Vector{0}}

\newcommand{\vs}{\Vector{s}}

\newcommand{\vw}{\Vector{w}}

\newcommand{\gP}{\Vector{\psi}}
\newcommand{\gPi}{\gP_i}
\newcommand{\gPii}{\gP_{i+1}}
\newcommand{\gPki}{\gP_{k,i}}
\newcommand{\gPkii}{\gP_{k,i+1}}
\newcommand{\gPs}{\gP^\star}

\newcommand{\gL}{\Vector{\lambda}}
\newcommand{\gLi}{\gL_i}
\newcommand{\gLki}{\gL_{k,i}}

\newcommand{\mI}{\Matrix{I}}

\newcommand{\mD}{\Matrix{D}}

\newcommand{\mP}{\Matrix{P}}
\newcommand{\mPi}{\mP_i}
\newcommand{\mPki}{\mP_{k,i}}
\newcommand{\mPb}{\overline{\mP}}

\newcommand{\mPt}{\Tilde{\mP}}
\newcommand{\mPtki}{\Tilde{\mP}_{k,i}}

\newcommand{\mG}{\Matrix{G}}

\newcommand{\mGb}{\overline{\mG}}

\newcommand{\mGbki}{\mGb_{k,i}}

\newcommand{\mA}{\Matrix{A}}

\newcommand{\opT}{\Operator{T}}

\newcommand{\opTki}{\opT_{k, i}}
\newcommand{\opTqi}{\opT_{i}^{(\sQ)}}
\newcommand{\opTqki}{\opT_{k, i}^{(\sQ_k)}}

\newcommand{\opTai}{\opT_{i}^{(\Theta)}}

\newcommand{\opTaki}{\opT_{k, i}^{(\Theta_k)}}

\newcommand{\opP}{\Operator{P}}

\newcommand{\opL}{\Operator{K}}

\newcommand{\opLki}{\opL_{k,i}}
\newcommand{\opV}{\Operator{V}}

\newcommand{\opZ}{\Operator{\Phi}}
\newcommand{\opZki}{\opZ_{k,i}}

\newcommand{\dmin}{\delta_{\min}}
\newcommand{\dmax}{\delta_{\max}}

\newcommand{\iN}{{i\in\Natural}}

\newcommand{\kN}{{k\in\sN}}
\newcommand{\oinO}{{\omega\in\Omega}}

\newcommand{\sigmaF}{\mathcal{F}}
\newcommand{\pspace}{(\Omega, \sigmaF, \probMeas)}
\newcommand{\Pas}{{\probMeas\text{-a.s.}}}

\newcommand{\HS}{\Real^M}

\newcommand{\VS}{\Real^{MN}}
\newcommand{\VM}{\Real^{MN \times MN}}

\newcommand{\rmu}{\mu}
\newcommand{\rmuki}{\rmu_{k,i}}
\newcommand{\pP}{\mathbb{P}}
\newcommand{\pPi}{\pP_i}

\newcommand{\hf}{\hat{f}}

\newcommand{\ey}{\est{y}}
\newcommand{\sQki}{\sQ_{k,i}}
\newcommand{\sQs}{\sQ^\star}
\newcommand{\sQt}{\Tilde{\sQ}}

\newcommand{\mGbkim}{[\mGb_i]_{(f_{km},:)}}
\newcommand{\mGbkimp}{[\mGb_i]_{(f_{km},f_{pq})}}
\newcommand{\mGbkimpd}{[\mGb_i]_{(f_{km},f_{p'q'})}}

\newcommand{\Expi}[1]{\Expp{i}{#1}}
\newcommand{\piki}{\pi_{k,i}}
\newcommand{\pikii}{\pi_{k,i+1}}
\newcommand{\tpi}{\pi}

\newcommand{\tpiki}{\tpi_{k,i}}

\newcommand{\pv}{\vartheta}
\newcommand{\pvv}{\Vector{\pv}}
\newcommand{\fq}{\Functional{W}}

\newcommand{\lSN}{\acute{\Natural}}
\newcommand{\lAN}{\Natural}

\newcommand{\iNt}{{i\in\lAN}}

\title{Distributed Fixed-Point Algorithms for Dynamic Convex Optimization over Decentralized and Unbalanced Wireless Networks%
\thanks{The authors acknowledge the financial support by the Federal Ministry of Education and Research of Germany in the program of ``Souver\"an. Digital. Vernetzt.'' Joint project 6G-RIC, project identification numbers: 16KISK020K and 16KISK030, and the German Research Foundation (DFG) within their priority program SPP 1914 ``Cyber-Physical Networking''.}%
}
\author{%
    Navneet Agrawal$^1$, Renato L.~G.~Cavalcante$^2$, S{\l}awomir Sta\'nczak$^{1, 2}$ \\
    $^1$\textit{Technische Universit\"at Berlin}, $^2$\textit{Fraunhofer Heinrich-Hertz-Institute}
}

\IEEEpubid{%
    \begin{minipage}{\columnwidth}
    \copyright 2024 IEEE. Personal use of this material is permitted. Permission from IEEE must be obtained for all other uses, in any current or future media, including reprinting / republishing this material for advertising or promotional purposes, creating new collective works, for resale or redistribution to servers or lists, or reuse of any copyrighted component of this work in other works.
    \end{minipage}%
    \hspace{1em}%
    \begin{minipage}{\columnwidth}
        $ $
    \end{minipage}
}

\begin{document}

% \address{$^1$Technische Universit\"at Berlin, $^2$Fraunhofer Heinrich-Hertz-Institute}
\maketitle

\IEEEpubidadjcol

\begin{abstract}
    We consider problems where agents in a network seek a common quantity, measured independently and periodically by each agent through a local time-varying process.
    Numerous solvers addressing such problems have been developed in the past, featuring various adaptations of the local processing and the consensus step.
    However, existing solvers still lack support for advanced techniques, such as superiorization and over-the-air function computation (OTA-C).
    To address this limitation, we introduce a comprehensive framework for the analysis of distributed algorithms by characterizing them using the quasi-Fej\'er type algorithms and an extensive communication model.
    Under weak assumptions, we prove almost sure convergence of the algorithm to a common estimate for all agents.
    Moreover, we develop a specific class of algorithms within this framework to tackle distributed optimization problems with time-varying objectives, and prove its convergence to a solution.
    We also present a novel OTA-C protocol for consensus step in large decentralized networks, reducing communication overhead and enhancing network autonomy as compared to the existing protocols.
    The effectiveness of the algorithm, featuring superiorization and OTA-C, is demonstrated in a real-world application of distributed supervised learning over time-varying wireless networks, highlighting its low-latency and energy-efficiency compared to standard approaches.
    
    % Typically, in solvers addressing problems of this type, each agent a local processing step followed by a consensus step.
    % and numerous algorithms have been developed in the past featuring various adaptations of these two steps.
    % We explicitly characterize the convergence point of our algorithm when employing a specific quasi-Fej\'er-type fixed-point algorithm, namely, the superiorized adaptive projected subgradient method (sAPSM), for a general class of dynamic convex optimization problems.
    % Numerous algorithms have been developed in the past, featuring various adaptations of these two steps for different applications.
    % This paper addresses dynamic distributed optimization problems over decentralized and time-varying networks with random directed graphs. 
    % The objective is for agents to reach an agreement on a common quantity, each agent independently measuring it through a local time-varying process. 
\end{abstract}


\section{Introduction} \label{sec:intro}
We consider distributed optimization problems in multiagent systems, where the underlying network is decentralized and time-varying, with possibly random and non-symmetric (directed) graph topologies.
The objective is for agents to reach agreement on the estimate of a common quantity, whose measurements are acquired by each agent independently and periodically via some local time-varying process.
Problems of this type have wide-ranging applications, including, adaptive control \cite{Nedic2018,Yang2020} and distributed learning \cite{Cavalcante2013}.

In general, solvers addressing such problems follow a two-step iterative approach: the \emph{local processing step}, where agents independently process their acquired measurements, followed by the \emph{consensus step}, where they communicate over the network to seek agreement.
Over the past four decades, numerous algorithms have been developed featuring various adaptations of these steps for different applications and/or system requirements \cite{Tsitsiklis1986,Nedic2009,Cavalcante2013,Nedic2014,Nedic2016,Li2020} (also see \cite{Nedic2018,Yang2019} and references therein).
However, some advanced techniques in machine learning and wireless communication that are known to be better suited for many current and envisioned applications are still not supported by the aforementioned studies.
Among these, two notable techniques are: \emph{Superiorization},\footnotemark~an efficient method to construct heuristics for constrained optimization problems \cite{Censor2010,Fink2022,Fink2023}; 
\footnotetext{\emph{Superiorization}~\cite{Censor2010} is a technique where a sequence of specifically designed perturbations are added to an iterative algorithm to steer its iterates towards vectors that are ``superior'' in terms of some desirable properties. A \emph{superiorizable} or bounded perturbation resilient algorithm enjoys convergence guarantees of the original iterative algorithm.}%
and \emph{over-the-air function computation} (OTA-C),\footnotemark~a scalable solution for distributed function computation over wireless networks \cite{Goldenbaum2013a,Agrawal2019}.
\footnotetext{In OTA-C~\cite{Goldenbaum2013a}, the additive structure of the wireless multiple access channels (WMAC) is exploited to compute certain functions, like consensus protocol, efficiently over large and dense wireless networks.}%
In \cite{Agrawal2023a}, we propose a class of distributed algorithms based on the adaptive projected subgradient method (APSM) \cite{Yamada2005} supporting both of these technologies.
However, their application is currently limited to the networks with symmetric (undirected) graph topologies.

In this paper, we overcome these limitations by making three key contributions.
% Therefore, in this paper, we extend the algorithm proposed in \cite{Agrawal2023a} to accommodate a broader range of optimization techniques and communication protocols.
First, we introduce a comprehensive framework of distributed algorithms that accommodates a broader range of optimization techniques and communication protocols.
We achieve this by characterizing the local processing step using operators generating quasi-Fej\'er monotone sequences, described later in Definition \ref{def:T_QF}.
%, which includes the superiorized APSM (sAPSM) \cite{Fink2023} as a special case (see Definition \ref{def:T_APSM}).
In the consensus step, we employ an abstract communication model that extends \cite{Agrawal2023a} to support random directed graphs.
% This characterization includes the superiorized APSM (sAPSM) \cite{Fink2023} used in \cite{Agrawal2023a} as a special case.
% , and provides a unified framework for the analysis of various existing algorithms \cite{Liu2017,Li2020}.
% Moreover, we extend the communication model of \cite{Agrawal2023a} to support a wider range of communication system, including the random directed graphs (see Assumption \ref{as:net}).
Second, under certain practical assumptions, we provide guarantees for convergence of the proposed algorithms in Theorem \ref{thm:T_QF}.
Our proofs
% \footnote{Complete proofs are provided at \texttt{https://t.ly/C9idY}.}~%
% (available in the arxiv version of the paper \cite{agrawal2024unbalanced_arxiv}) 
utilize \emph{time-varying quadratic Lyapunov functions} \cite{Touri2013,Nedic2016}, making the analysis considerably more involved than the undirected case in \cite{Agrawal2023a}.
% In Theorem \ref{thm:T_QF}, we establish almost sure convergence of the proposed scheme \eqref{eq:agent_scheme} using operators generating quasi-Fej\'er monotone sequences.
We also present a specific algorithm based on the superiorized APSM (sAPSM) to tackle a class of dynamic convex optimization problems, and, in Theorem \ref{thm:T_APSM}, prove its convergence to a unique solution.
Third, in Section \ref{sec:comm}, we propose a novel OTA-C protocol that reduces the communication overhead and grants more autonomy to the agents, compared to previous studies \cite{Nedic2009,Cavalcante2013,Agrawal2023} that require some prior coordination and additional overhead to ensure that the underlying graph remains balanced or undirected.
In Section \ref{sec:sim}, the proposed algorithm, featuring sAPSM and OTA-C, is applied to a distributed supervised learning problem, and shown to outperform the standard approaches.

%\vspace{-1em}
\emph{Notation}: %\label{sec:prelims}
The space of complex, real, nonnegative real, and natural numbers (including zero) are given by $\Complex$, $\Real$, $\Realp$, and $\Natural$, respectively.
Symbols $\vOne_M$, $\mI_M$ and $\otimes$ denote the vector of $M$ ones, the $M\times M$ identity matrix, and the Kronecker product, respectively.
Any $M$-dimensional real-valued vector belongs to the Hilbert space $(\Real^M, \innerprod{\cdot}{\cdot})$ with inner-product $(\forall \vv,\vy \in \Real^M)$ $\innerprod{\vv}{\vy} := \vv^T\vy$ and induced norm $\norm{\vv} := \sqrt{\vv^T \vv}$.
We denote by $\ell_1^+$ the space of all nonnegative sequences of real numbers, such that for any sequence $(\xi_i)_\iN\in\ell_1^+$ we have $\sum_\iN \xi_i < \infty$.
The underlying probability space is $\pspace$,
%where $\Omega$ is an arbitrary sample set, $\mathcal{F}$ is a $\sigma$ algebra of subsets of $\Omega$, and $\probMeas$ is a probability measure.
and a \emph{random variable} is a measurable map from $\Omega$ to some real vector space.
% Expressions involving random variables must be understood to hold $\mathbb{P}$-almost surely ($\Pas$).
Given random variables $\vx$ and $\vy$, the conditional expectation of $\vy$ w.r.t.~the sigma algebra generated by $\vx$ is denoted by $\Exp{\vy \mid \vx}$.
% A matrix $\mA$ is called \emph{stochastic} if its elements are nonnegative, and its rows sum up to one.
% A directed graph $\sG=(\sN, \sE)$ is called \emph{weakly} (resp., strongly) connected if it is possible to reach any node starting from any other node by traversing the graph over undirected (resp., directed) edges in $\sE$.

% \clearpage

\section{General class of distributed algorithms}\label{sec:general}
%\vspace{-0.7em}
In this paper, we study the general class of distributed algorithms that follows an \emph{adapt-then-combine} \cite{Sayed2014} strategy.
% In this paper, we strive to unify these approaches to facilitate development and analysis of 
Solvers of this type employ a two-step iterative approach: a local optimization step, often implementing a fixed-point algorithm \cite{Bertsekas1983, Tsitsiklis1986, Fullmer2016, Liu2017, Li2020}, followed by a consensus step that fosters agreement among agents using networking protocols such as broadcast or gossip protocols in wireless networks \cite{Nedic2009, Nedic2010, Cavalcante2013, Nedic2018, Agrawal2023}.
In this section, we formulate these two steps to accommodate diverse optimization and communication techniques, establishing a unified framework for representing and analyzing a broad class of distributed optimization algorithms, including those employed in the aforementioned studies.
Later in Section \ref{sec:APSM}, we also develop a particular algorithm under this framework that addresses a class of time-varying distributed optimization problems, and strengthen the convergence results presented in this section.

%\vspace{-1.0em}
\subsection{System description and agent dynamics}\label{sec:prob}
%\vspace{-0.5em}
We consider a multiagent system consisting of a network of $N$ agents, % where agent dynamics is given by the following two-step iterative scheme.
% represented by a time-varying directed graph $\sGi := (\sN, \sEi)$, where $\iN$ is the time-index, $\sN := \{1,\dots,N\}$ is the set of graph vertices, and $\sEi \subseteq \sN\times\sN$ is the set of directed edges.
where each agent $\kN:=\{1, \dots, N\}$ implements the following two-step dynamics at every time step $\iN$:
%\vspace{-0.25em}
\begin{subequations}\label{eq:agent_scheme}
    \begin{align}
        \gLki               &=    \opP_{\sX}(\opTki(\gPki)) \label{eq:agent_scheme_local_opt} \\
        \gP_{k,i+1}         &=    (1 - \beta_i)\ \gLki + \beta_i\ \opLki(\gL_i), \label{eq:agent_scheme_consensus}
    \end{align}
\end{subequations}
%\vspace{-1.5em}
where, starting with an arbitrary initialization $\gP_{k,0}$, the current estimate of agent $k$ at time $i$ is $\gPki\in\HS$.
The vector $\gL_i$ is formed by stacking $\gLki$ for all agents $\kN$, and the sequence $(\beta_i)_\iN$ is a design parameter defined as $\beta_i = (i+1)^{-\alpha}$ with ${0.5 < \alpha \leq 1}$.
Note that $(\beta_i)_\iN \notin \ell_1^+$ as the series $\sum_\iN \beta_i$ diverges, but $(\beta_i^2)_\iN \in \ell_1^+$.
In the \emph{local processing step} \eqref{eq:agent_scheme_local_opt}, the current estimate $\gPki$ is first updated by applying the operator $\opTki:\HS\to\HS$ that embeds local information only available to agent $k$ at time $i$ via, for instance, some measurement process.
It is followed by the mapping $\opP_{\sX}$ which projects\footnotemark~any vector in $\HS$ onto a closed and convex set $\sX\subset\HS$
The mapping $\opP_{\sX}$ enforces a technical requirement for the exchange of $\gLki$ over a finite-capacity communication system.
We characterize the sequence of operators $(\opTki)_\iN$ explicitly later in Section \ref{sec:T_char}. % as a sequence of operators generating quasi-Fej\'er monotone sequences
\footnotetext{In Hilbert space $\Space{H}$, projection of $\vy\in\Space{H}$ onto a closed convex set $\sX\subset\Space{H}$ is the unique $\opP_{\sX}(\vy)\in\sX$ satisfying $\min\limits_{\vx\in\sX} \norm{\vy - \vx}_{\Space{H}} = \norm{\vy - \opP_{\sX}(\vy)}_{\Space{H}}$.}

In the \emph{consensus step} \eqref{eq:agent_scheme_consensus}, agent $k$ incorporates information received from its neighbors over the network.
The information exchange over the network is modeled by the random mapping $\opLki$ which, for almost every (a.e.) $\oinO$, takes the form: $(\forall \iN) (\forall \kN) (\gP\in\VS)$
% In the following, we define $\opLi$ as the random map formed by stacking $(\opLki)_\kN$ for all agents.
% More precisely, for any event $\oinO$ and input $\gP\in\VS$, we define $\opLi(\omega, \gP) := [(\opL_{1,i}(\omega, \gP))^T, \dots, (\opL_{N,i}(\omega, \gP))^T]^T \in \VS$.
% The mapping $\opLi$ thus formed, models the communication between all agents in the graph $\sGi$.
% %\vspace{-0.25em}
\begin{align} \label{eq:comm_model}
   \opLki(\omega, \gP) = \mPki(\omega)\ \gP + \vnki(\omega), %\quad \text{ for a.e. } \oinO,
\end{align}
%\vspace{-1.25em}
where $\mPki$ is a random matrix taking values in $\Real^{M\times MN}$, and $\vnki$ is a random vector in $\Real^M$ that may depend on the input $\gP$.
% In other words, $\opLki$ in \eqref{eq:comm_model} linearly transforms the input $\gL$ by randomly weighing its every component with corresponding edge-weight on the graph, given by the respective element of $\mPki$, while introducing an additive state-dependent noise $\vnki$.
Note that the model in \eqref{eq:comm_model} covers a large class of wireless protocols with digital and analog transmissions, including the novel OTA-C protocol proposed in Section \ref{sec:comm}.
% Let $\mPi$ and $\vni$ be stacked versions of $\mPki$ and $\vnki$, respectively, for all $\kN$.

The network is represented by the directed graph $\sGi:=(\sN, \sE_i)$ at any time $\iN$, where $\sE_i\in\sN\times\sN$ is the set of directed edges between the agents.
Let $\mPi\in\VM$ be the matrix formed by stacking $\mPki$ for all $\kN$.
Assuming that each coordinate of vector $\gLki\in\HS$ is exchanged over a separate i.i.d.~realization of the random graph $\sGi$,\footnotemark~the expectation of the random matrix $\mPi$ takes the form $\mPb_i := \Exp{\mPi} = \mA_i \otimes \mI_M$, where $\mA_i\in\Realp^{N\times N}$ and, for any $(p,q)\in\sN\times\sN$, the scalar $[\mA_i]_{(p,q)}$ represents the expected edge-weight of the edge from agent $q$ to $p$ in the graph $\sGi$.

\footnotetext{In a random graph, the weights over its edges are random variables.}

We make the following assumptions on the sequence of graphs $(\sG_i)_\iN$ and corresponding matrices $(\mA_i)_\iN$.
%\vspace{-0.5em}
\begin{assumption}[Network assumptions] \label{as:net} For all $\iN$:
    % We assume that the sequence of matrices $(\mA_i)_\iN$ satisfy the following conditions:

    % The following holds for all $\iN$:

    (i) The matrix $\mA_i$ is stochastic,\footnote{A matrix $\mA\in\Real^{N\times N}$ is called \emph{stochastic} if $\mA\geq 0$, and $\mA\vOne_N=\vOne_N$.}~and it is compliant with the graph $\sGi$, i.e., $[\mA_i]_{(k,l)} > 0$ if and only if $(l,k)\in\sE_i$.

    (ii) There exists $\epsilon>0$ such that $[\mA_i]_{(k,k)} \geq \epsilon$ for all $\kN$, and $[\mA_i]_{(k,l)} \geq \epsilon$ for all $(l,k)\in\sE_i$.

    (iii) The graph $\sGi$ is \emph{strongly connected} in expectation, i.e., there exist $n>0$, $n\in\Natural$, such that $(\mA_i)^n$ is a positive matrix.
    % \footnote{A directed graph $\sGi$ with random edge-weights is said to be \emph{strongly connected in expectation} if, in the graph with expected edge-weights, every pair of agents is connected via a sequence of directed edges.}

    (iv) $\Exp{\vni \mid \gPi} = \vzero$, $\Pas$
    
    (v) ${\Exp{\norm{\mPi^T\mPi}_2} < \infty}$, ${\Exp{\norm{\vni}^2 \mid \gPi} < \infty}$, and ${\Exp{\norm{\mPi^T \vni} \mid \gPi} < \infty}$, $\Pas$
    
\end{assumption}
%\vspace{-0.5em}

% \footnotetext{A matrix is \emph{stochastic} if it is nonnegative, and all its rows sum up to one.}
% \footnotetext{A directed graph with random edge-weights is said to be \emph{strongly connected in expectation} if, in the graph with expected edge-weights, every pair of agents is connected via a sequence of directed edges.}

\begin{remark} \label{rem:model_assumptions}
    Assumptions \ref{as:net}(i) and \ref{as:net}(iii) are standard in literature \cite{Nedic2014,Nedic2016,Li2020} and \ref{as:net}(ii) facilitates analysis of systems over directed graphs \cite{Touri2013,Nedic2016}.
    % According to \cite{Wu2005a,Nedic2016}, for the sequence of matrices $(\mA_i)$ under Assumption \ref{as:net}, there exists a sequence of vectors $(\pi_i)$ in $\Real^N$ such that $\pi_i^T = \pi_{i+1}^T \mA_i$.
    % This property is essential for our convergence proofs, which utilize a time-varying quadratic Lyapunov function that depends on $\pi_i$ at any $\iN$.
    % The statistical conditions on the communication model in Assumption \ref{as:net}(iv) are satisfied, in general, for wireless network protocols.
    Moreover, Assumptions \ref{as:net}(iv) and \ref{as:net}(v) are weaker than their counterparts in \cite{Agrawal2023a}, and hence, the model in \eqref{eq:comm_model} under Assumption \ref{as:net} is more general.
\end{remark}

%\vspace{-1.2em}
\subsection{Quasi-Fej\'erian characterization of ${(\opTki)_\iN}$ and analysis} \label{sec:T_char}
%\vspace{-0.5em}
Most fixed-point algorithms (e.g., those with nonexpansive operators) and their variants (e.g., superiorized APSM \cite{Fink2023}) generate quasi-Fej\'er monotone sequences (QFMS) of type-III \cite{Combettes2001,Bauschke2017}.
Hence, we characterize the sequence of operators $(\opTki)_\iN$, for every $\kN$, as generators of QFMS in the sense of Definition \ref{def:T_QF} below.
In Theorem \ref{thm:T_QF}, we establish sufficient conditions for almost sure convergence of the algorithm based on these operators to a common point for all agents.
% Subsequently, we particularize these results for a variant of $\opTki$ based on the sAPSM operator $\opTaki$, and prove in Theorem \ref{thm:T_APSM} that the convergence point is indeed a solution of the corresponding dynamic convex optimization problem.
% Sufficient conditions for almost sure convergence are established for both cases in Theorem \ref{thm:T_QF} and Theorem \ref{thm:T_APSM}, respectively.

%\vspace{-0.5em}
\begin{definition}[$(\opTqi)_\iN$ : Quasi-Fej\'er monotone sequence (QFMS) generator] \label{def:T_QF}
    The sequence of operators $(\opTqi)_\iN$ is called a \emph{QFMS generator} w.r.t.~a nonempty set $\sQ\subset\HS$ if,
    for the sequence $(\vx_i)_\iN$ of vectors generated via $(\forall\iN)\ \vx_{i+1} = \opTqi(\vx_i)$, $\vx_0\in\HS$, and for any $\vx\in\sQ$, there exists $(\epsilon_{i})_\iN \in \ell_1^+$ such that:
    %\vspace{-0.3em}
    \begin{align} \label{eq:quasiF}
        (\forall\iN) \qquad \norm{\vx_{i+1} - \vx}^2 \leq \norm{\vx_i - \vx}^2 + \epsilon_{i}.
    \end{align}
    % \vspace{-1.3em}
    
    % \noindent Then, the sequence of operators $(\opTqi)_\iN$ is called a \emph{quasi-Fej\'er monotone sequence generator}.
    % Starting from an arbitrary $\vx_0\in\HS$, the sequence $(\vx_i)_\iN$ generated by $(\forall \iN)\ \vx_{i+1} = \opTqki(\vx_i)$ is a \emph{quasi-Fej\'er monotone sequence of type-III} \cite{Combettes2001,Fink2023} with respect to some nonempty set $\sQ_k \subset \sX$.
    % More precisely, for any $\vx\in\sQ_k$, there exists $(\epsilon_{k,i})_\iN \in \ell_1^+$ such that:

    % Note that, for all $\iN$, $\opTki(\vx_i)$ and $\vx$ are random vectors in $\Real^M$, and $\epsilon_{k,i}$ is a random variable in $\Realp$.
    % The mapping corresponding to an agent $n\in\sN$ is denoted by $\opT_{i,n}$, and the respective set by $\sQ_n$.
    % Furthermore, agents have the knowledge of an orthogonal projection operator $\opP_{\sX}:\Real^M\to\Real^M$ onto a nonempty closed convex set $\sX\subset \Real^M$.
    % \noindent In addition, the set $\sQs := \bigcap_{\kN} \sQ_k \cap \sX$ is nonempty.
\end{definition}

% The following result establishes sufficient conditions for convergence of scheme \eqref{eq:agent_scheme} to a common point in $\sX$ for every agent $\kN$, where each agent $\kN$ uses the sequence of operators $(\opTqki)_\iN$ as in Definition \ref{def:T_QF} for the local optimization step \eqref{eq:agent_scheme_local_opt}.

% Theorem \ref{thm:T_QF} below establishes sufficient conditions for convergence of the scheme \eqref{eq:agent_scheme} to the same point for all agents, where each agent $\kN$ employs a QFMS generator $(\opTqki)_\iN$ in step \eqref{eq:agent_scheme_local_opt}.

\begin{theorem} \label{thm:T_QF}
    Suppose that Assumption \ref{as:net} holds in a system where each agent $\kN$ implements the scheme in \eqref{eq:agent_scheme}, where $\opLki$ is given by \eqref{eq:comm_model}, and, for all $\kN$, $(\opTki)_\iN \equiv (\opTqki)_\iN$ is a QFMS generator w.r.t.~set $\sQ_k\subset\HS$ as defined in Definition \ref{def:T_QF}.
    Moreover, assume that the set $\sQs := \bigcap_{\kN} \sQ_k \cap \sX \subset \HS$ is nonempty.
    Then, each of the following statements hold:
    
    (i) \emph{(Convergence)}:
        % Let $\gPi\in\VS$ be the vector formed by stacking $\gPki\in\HS$ for all $\kN$.
        For any $\gPs \in \sQs$, the sequence $(\norm{\gPki - \gPs}^2)_\iNt$ converges ($\Pas$) for all $\kN$.
        % where $(\tpi_i)_\iN$ is a unique sequence of stochastic vectors such that ${(\forall \iN)(\forall \kN)(\exists \delta>0)\ \tpiki \geq \delta}$ and $(\forall \iN)$ ${\sum_\kN \tpiki = 1}$.
        Hence, every $(\gPki)_\iNt$ is bounded $\Pas$, and has an accumulation point.
            
    (ii) \emph{(Consensus)}: All agents in $\sN$ reach consensus, i.e.
    %\vspace{-0.5em}
    \begin{align*}
        (\forall (p,q) \in\sN\times\sN)(\iNt) \quad \lim_{i\to\infty} \norm{\gP_{p,i} - \gP_{q,i}} = 0,\quad \Pas
    \end{align*}
    %\vspace{-1.5em}

    (iii) \emph{(Characterization of accumulation points)}:
        In addition, assume that 
        the set $\sQs$ has a nonempty interior, i.e., for some $\vut\in\sQs$, $\exists \varrho > 0$ such that ${\sQs \supset \{\vu\in\Real^M \mid \norm{\vu - \vut} \leq \varrho\} \neq \emptyset}$.
        
        %the following statements hold:
        
        % (a) Set $\sQs$ has a nonempty interior, i.e., for some $\vut\in\sQs$, $\exists \varrho > 0$ such that ${\sQs \supset \{\vu\in\Real^M \mid \norm{\vu - \vut} \leq \varrho\} \neq \emptyset}$, and
        
        % (b) Sequence $(\sum_\kN \norm{\gPki - \gPs}^2)_\iN$ is convergent $\Pas$
        
        Then, for all agents $\kN$, the sequence $(\gPki)_\iNt$ converges to the same point in $\sX$, $\Pas$

    % \begin{proof}
    %     \input{proof_thm1.tex}
    % \end{proof}
    \begin{proof}
        Proof given in the Appendix \ref{thm:1}.
        % Proof is provided in the extended version \cite{Agrawal2024distributed}.
    \end{proof}
\end{theorem}

\begin{remark}\label{rem:T_Q}
    The condition in \eqref{eq:quasiF}, known as quasi-Fej\'er monotonicity (QFM), holds for several fixed-point algorithms, and it has proven to be an efficient tool for their analysis \cite{Combettes2001,Combettes2015}.
    However, the analysis of distributed algorithms based on the QMF condition is a novel contribution of this paper.
    Note that Theorem \ref{thm:T_QF} falls short of an explicit characterization of the point of convergence, for instance, to the set $\sQs$.
    In Section \ref{sec:APSM}, we develop a variant of the QFMS generator, for which the point of convergence of the scheme \eqref{eq:agent_scheme} can be characterized as a time-invariant solution of an infinite sequence of time-varying convex optimization problems.
\end{remark}

% To provide a concrete problem and algorithm, we consider the \emph{adaptive projected subgradient method} (APSM), with a sequence of bounded perturbation introduced to its iterates \cite{yamada2005adaptive,Fink2023fixed}.
% Note that APSM is an $\eta$-attracting quasi-nonexpansive mapping \cite[Proposition 2]{yamada2005adaptive}.

% \section{Solving distributed optimization problems}\label{sec:appl}

%\vspace{-1em}
\subsection{Dynamic distributed convex optimization via sAPSM}\label{sec:APSM}
%\vspace{-0.5em}
In this section, we develop algorithms to solve a class of distributed convex optimization problems with time-varying objectives.
The proposed algorithm is a variation of \eqref{eq:agent_scheme} with a specific QFMS generator based on the superiorized APSM (sAPSM) (see Definition \ref{def:T_APSM}).
% In the following, 
% Then, we define the sAPSM, which, when employed in the scheme \eqref{eq:agent_scheme}, allows agents to reach a common solution of the problem.
We consider the following problem $\pPi$ at any time $\iN$:
%\vspace{-0.3em}
\begin{equation}
    \label{eq:optprob}
    \begin{aligned} 
        \pPi: \quad \underset{\gP_1\in\sX, \dots, \gP_N\in\sX}{\text{minimize}} \sum_{k\in\sN} \cT(\gP_k), \quad
        \text{s.t.}\,\,\, \gP_1=\dots=\gP_N,
    \end{aligned}
\end{equation}
where, for all $\kN$, the cost function $\cT:\HS\to\Realp$ is convex and possibly nonsmooth with $\min_{\vx\in\sX} \cT(\vx) = 0$.\footnote{This is a common assumption in subgradient literature \cite{Yamada2005,Nedic2009,Cavalcante2013,Fink2023}.}~%
Define $\sQ^{(i)} := \bigcap_{k\in\sN} \sQki$, $\sQki := \{\vh\in\sX \mid \cT(\vh) = 0\}$, and assume that $\sQs := \bigcap_{\iN} \sQ^{(i)} \neq \emptyset$.
Ideally, the objective of agents is to find a solution to all infinitely many problems $(\pPi)_\iN$, assuming that such a solution exists.
However, due to system causality and limited memory, finding such a point is practically infeasible.
Instead, we relax the problem to finding a point in the set of solutions to all but finitely many problems $(\pPi)_\iN$, that is \cite{Yamada2005,Cavalcante2013,Agrawal2023a}:
%\vspace{-0.5em}
\begin{align}
    \label{eq:relax_set}
    \text{Find}\ \bh \in \sQt := \overline{\liminf_{i\to\infty} \sQ^{(i)}} = \overline{\bigcup_{n\in\Natural}\bigcap_{i\geq n} \sQ^{(i)}} \supset \sQs \neq \emptyset,
\end{align}
where $\overline{\sC}$ denotes the closure of the set $\sC$.

The APSM generates a sequence of estimates that are known to converge in $\sQt$ \cite{Yamada2005,Cavalcante2013}.
Moreover, the APSM is superiorizable, i.e., it is resilient to bounded perturbations.
In the following, we first define the sAPSM operators in Definition \ref{def:T_APSM} below, and then, in Theorem \ref{thm:T_APSM}, prove that the sAPSM based scheme \eqref{eq:agent_scheme} solves \eqref{eq:relax_set}.
Note that the sequence $(\opTaki)_\iN$ generated by the sAPSM operators is a QFMS generator w.r.t.~set $\sQ_k = \cap_\iN \sQki$ \cite{Fink2023}.

%\vspace{-0.5em}
\begin{definition}[$(\opTai)_\iN$ : Superiorized APSM (sAPSM) sequence generator]\label{def:T_APSM}
    Given a sequence of convex cost function $(\Theta_i)_\iN$, where each $\Theta_i:\HS\to\Realp$ and $\min_{\vx\in\sX} \Theta_i(\vx) = 0$, the sequence of mappings $(\opTai)_\iN$ generating $(\vx_i)_\iN$ via: %$(\forall \iN)\ \vx_{i+1} = \opTai(\vx)$, where each $\opTai$ is defined as:
    %\vspace{-0.5em}
    \begin{align}\label{eq:T_APSM}
        (\forall \iN) \qquad \vx_{i+1} = \opTai(\vx_i) := \vx_i - \opZ_i(\vx_i) + \zeta_i\vz_i,
    \end{align}
    is called a \emph{sAPSM sequence generator}, where
    $(\zeta_i \vz_i)_\iN$ is a sequence of bounded perturbations\footnotemark~in $\HS$, and
    \footnotetext{A sequence $(\zeta_i \vz_i)_\iN$ in $\HS$ is called a \emph{sequence of bounded perturbations} if $(\zeta_i)_\iN \in \ell_1^+$, and $(\exists r>0)(\forall\iN)\ \norm{\vz_i}\leq r$.}%
    $\opZ_i:\HS\to\HS$ is defined as: $(\forall \iN)(\vx\in\HS)$
    %\vspace{-0.5em}
    \begin{align}\label{eq:subGop}
        \opZ_i(\vx) := (\rmu_i\Theta_i(\vx)/\norm{\Theta'_i(\vx)}^2)\  \Theta'_i(\vx),
    \end{align}
    if $\norm{\Theta'_i(\vx)} \neq 0$, otherwise $\opZ_i(\vx) = 0$,
    where $\rmu_i \in (0,2)$ is a design parameter, and
    $\Theta'_i(\vx) \in\partial \Theta_i(\vx)$.\footnotemark
\end{definition}
\footnotetext{Given a convex function $\Theta:\HS\to\Real$ at $\vx\in\HS$, we define $\partial\Theta(\vx) := \{\vh\in\HS \mid (\vy - \vx)^T \vh + \Theta(\vx) \leq \Theta(\vy), \, \, \forall\vy\in\HS \}$.}

% The following assumption is standard in the literature related to the subgradient methods \cite{Yamada2005,Nedic2009,Cavalcante2013}.
\noindent We further make two standard technical assumptions \cite{Yamada2005,Cavalcante2013}.
%\vspace{-0.5em}
\begin{assumption}[Problem-specific assumptions]\label{as:bound_sugGop}
    A time-invariant solution to all problems in $(\pPi)_\iN$ exists, i.e., $\sQs \neq \emptyset$, and, for all $\kN$, $(\cS(\gPki))_\iN$ is bounded $\Pas$
\end{assumption}

% The following result establishes, in particular, sufficient conditions for all agents to asymptotically reach consensus on a solution to \eqref{eq:relax_set}.

\begin{theorem}\label{thm:T_APSM}
    Suppose that Assumptions \ref{as:net} and \ref{as:bound_sugGop} hold in a system where, given $\Theta_k := (\cT)_\iN$, each agent $\kN$ implements the scheme in \eqref{eq:agent_scheme} where $\opLki$ is given by \eqref{eq:comm_model} and, for all $\kN$, $(\opTki)_\iN \equiv (\opTaki)_\iN$ is a sAPSM sequence generator as defined in Definition \ref{def:T_APSM}.
    Then, in addition to the results already established in Theorem \ref{thm:T_QF}, the following statements hold:
    % for all $\kN$, the sequence $(\gPki)_\iN$ generated by \eqref{eq:agent_scheme} satisfies the following conditions, in addition to the results already established in Theorem \ref{thm:T_QF}:

    % (i) The sequence generated by $(\opTaki)_\iN$ is a QFMS generator w.r.t.~set $\sQ_k = \cap_\iN \sQki$.
    % Sufficient conditions for Theorem \ref{thm:T_QF} are satisfied. In particular, for all $\kN$, sAPSM $(\opTaki)_\iN$ is a QFMS generator.
    
    (i) The sequence $(\gPki)_\iN$ asymptotically minimize the local cost functions, i.e.,
        %\vspace{-0.5em}
        \begin{align*}
            (\forall \kN)\quad \lim\limits_{i\to\infty} \cT(\gPki) = 0, \quad \Pas
        \end{align*}
        %\vspace{-1.5em}

    (ii) For an interior point $\vut\in\sQs$, let us define sets $\sS_1 := \{\iN \mid \sum_{k\in\sN} \min_{\vx\in\lev_{\leq 0} \cT} \norm{\gPki-\vx} > \vartheta\}$,\footnotemark~and $\sS_2 := \{\iN \mid \sum_{k\in\sN} \norm{\vut - \gPki} \leq r\}$, $\Pas$~In addition to Assumptions \ref{as:net} and \ref{as:bound_sugGop}, suppose that: $(\forall \vartheta>0, \forall r>0, \exists\xi>0)$
    \footnotetext{Given a convex function $\Theta:\HS\to\Realp$ and $c\geq 0$, we define $\lev_{\leq c} \Theta := \{\vh\in\HS \mid \Theta(\vh) \leq c\}$.}
    % for all $\vartheta>0$, $r>0$, and for almost every $\omega \in \Omega$, there exists $\xi(\omega) >0$ such that
    %\vspace{-0.5em}
    \begin{align} \label{eq:complex_cond}
        \inf_{i\in\sS_1\cap\sS_2} \sum_{k\in\sN} \cT(\gPki) \geq \xi, \quad \Pas
    \end{align}
    %\vspace{-1.5em}
    
    % \noindent\sloppy where, for a.e.~$\oinO$, we define $\sS := \sS_1 \cap \sS_2$ with $\sS_1 := \{\iN \mid \sum_{k\in\sN} \min_{\vx\in\lev_{\leq 0} \cT} \norm{\gPki-\vx} > \vartheta\}$, $\sS_2 := \{\iN \mid \sum_{k\in\sN} \norm{\vut - \gPki} \leq r\}$, and, for all $\iN$ and $\kN$, $\lev_{\leq c} \cT := \{\vh\in\HS \mid \cT(\vh) \leq c\}$.
    % where $\sS(\omega) := \{\iN \mid \sum_{k\in\sN} \min_{\vx \in \sQki} \norm{\gPki(\omega) - \vx} > \vartheta \text{ and } \sum_{k\in\sN} \norm{\vut - \gPki(\omega)} \leq r\}$, 
    \noindent Then, for all $\kN$, the sequence of estimates $(\gPki)_\iN$ generated by \eqref{eq:agent_scheme} converges to a solution of \eqref{eq:relax_set} $\Pas$

    \begin{proof}
        Proof given in the Appendix \ref{thm:2}.
    \end{proof}
\end{theorem}

\begin{remark} \label{rem:T_APSM}
    Theorem \ref{thm:T_APSM} not only guarantees asymptotic convergence to a point that minimizes the cost $\cT$ for each agent $\kN$, but also characterizes the point of convergence explicitly as a solution of Problem \eqref{eq:relax_set}.
    % In simple words, condition \eqref{eq:complex_cond} states that the sequence $(\gPki)_\iN$ does not converge too quickly to the set $\sQt$.
    Note that the point of convergence in Theorem \ref{thm:T_APSM} is a $\sQt$-valued random variable.
    Hence, although the problem itself is deterministic, different runs of \eqref{eq:agent_scheme} may lead to different estimates in $\sQt$.
\end{remark}

%\vspace{-0.5em}
\section{OTA-C for decentralized consensus}\label{sec:comm}
%\vspace{-0.5em}
% In this section, we provide a brief description of the OTA-C based communication protocol for implementing the consensus step \eqref{eq:agent_scheme_consensus} in practical wireless systems.
In this section, we present a novel OTA-C protocol that extends our prior work \cite{Agrawal2023,Agrawal2023a} with some notable differences, as mentioned in Remark \ref{rem:OTAC}.
% Specifically, we relax the requirement in \cite[Def.3.1]{Agrawal2023a}, resulting in reduced communication overhead and fewer constraints on the WMAC.
% Also, in contrast to \cite{Agrawal2023a}, support for directed graphs allows agents to independently select their transmit powers. %, as long as it can be ensured that the resulting graph remains strongly connected in expectation (Assumption \ref{as:net}(iii)).
% previously in \cite[Def.3.1]{Agrawal2023a}, each row of the random matrix $\mPi$ is required to sum up to one, $\Pas$.
% Relaxing this requirement (see Assumption \ref{as:net}) results in an OTA-C protocol with reduced communication overhead and fewer constraints on the underlying WMAC.
% Furthermore, the adaptation to directed graphs gives more autonomy to the agents by allowing them to choose their transmit powers independently.
% This is in contrast to \cite{Agrawal2023a}, where equal transmit powers are required to ensure undirected graphs $(\sG_i)_\iN$ in expectation.
%The protocol design is based on the observation that, for executing the consensus step, an agent only needs to know a weighted combination of its neighbors' inputs, where the weights could be arbitrary and unknown \cite{kar2008distributed,molinari2022over,agrawal2023distributed}.
First, we provide a brief overview of the OTA-C protocol for implementing the consensus protocol \eqref{eq:agent_scheme_consensus} over a random directed graph (the protocol is given later in \eqref{eq:agents_protocol}).
% The extension to $M$-dimensional inputs, as discussed in Section \ref{sec:prob}, is achieved over $M$ i.i.d.~realization of the random graph $\sGi$ at every time $\iN$.
% % Note that the linear weights in the sum are unknown as they are functions of the random channel fading.
% The protocol is extended in Section \ref{sec:deploy_OTAC} to allow all agents in the network to obtain estimates for the vector-valued inputs, and this information is then used by the agents to implement the consensus protocol in \eqref{eq:scheme_consensus}.
% % Then, the consensus step for every agent in the network is combined to implement the consensus step in \eqref{eq:scheme_consensus}.
% % Note that \eqref{eq:ota_protocol} is equivalent to the consensus step in \eqref{eq:scheme_consensus} with $\check{\opL}_i$, instead of $\opLi$ in \eqref{eq:comm_model}, that models the information exchange in a network using the proposed OTA-C protocol.
Then, in Proposition \ref{prop:OTAC}, we prove that the sufficient conditions in Theorem \ref{thm:T_QF} and \ref{thm:T_APSM} are satisfied for the proposed OTA-C protocol based consensus.

Consider agent $r\in\sN$ and its (inward) neighbors in the set $\sN_r := \{k\in\sN \mid (k, r) \in \sE\}$.
% As discussed in Section \ref{sec:prob}, the $M$-dimensional inputs $\gLki$ are exchanged by implementing the protocol $M$ times independently at every iteration $\iN$.
As the protocol remains the same for every iteration $\iN$, we omit index $i$ in the following.
For every $m$th realization of the random graph, for $m\in\sM:=\{1,\dots,M\}$, any transmitting agent $k\in\sN_r$ generates a sequence of $B$ complex-valued random numbers $(s_k(1),\dots,s_k(B))$ as follows: 
for all $b=1,\dots,B$ and $\kN$, $s_k(b) := \sqrt{g_k(\lambda_k)} U_k(b)$, where $\lambda_k := \gLki^{(m)}$ denotes the $m$th element of $\gLki$, function $g_k(x) := P_k (x - \dmin)/(\dmax - \dmin)$, $(\dmax,\dmin):=(\max\sX, \min\sX)$, 
and $U_k$ is an i.i.d.~complex-valued random variable with $\abs{U_k(b)}=1$, and $\Exp{U_k} = 0$.
In addition, once every $\iN$, agents transmit $B'$ random numbers $(s'_k(1),\dots,s'_k(B'))$, encoded as above with a constant value $(\forall\kN)\ \lambda_k=\dmax$.
The signal received by any agent $r\in\sN$ over WMAC due to simultaneous transmissions by all $\sN_r$ is modeled as \cite{Goldenbaum2013a}: $(\forall b=1, \dots, B)$
%\vspace{-0.25em}
\begin{align} \label{eq:WMAC}
    q_r(b) = \sum_{k\in\sN_r} \xi_{kr}(b) s_k(b) + w_r(b),
\end{align}
%\vspace{-1.25em}
where $\xi_{kr}$ and $w_r$ are complex-valued random variables representing the fading channel between $k$ and $r$, and the receiver noise.
% We make the following assumptions on the WMAC model.
%\vspace{-0.5em}
\begin{assumption}[WMAC assumptions]\label{as:WMAC}
    For all agents $(r,t)\in\sN\times\sN$, and their (inward) neighbors $(k,l)\in \sN_r\times\sN_t$, the following properties hold:
    
    (i) \emph{Channel:} $\Exp{\abs{\xi_{kr}}^2} < \infty$, $\Exp{\abs{\xi_{kr}}^2 \abs{\xi_{lt}}^2} < \infty$.
    
    (ii) \emph{Noise:} $\Exp{w_r} = 0$, $\Exp{\abs{w_r}^2} < \infty$, $\Exp{\abs{w_r}^2 \abs{w_t}^2} < \infty$.
    
    (iii) The random variables $\xi_{rk}$ and $w_r$ in the WMAC model \eqref{eq:WMAC} are independent, and their statistics remains the same for the duration of $(MB + B')$ symbols.
\end{assumption}
%\vspace{-0.5em}

In response, at any iteration $\iN$ of the scheme, agent $r\in\sN$ receives $MB + B'$ symbols, where each received symbol takes the form given in \eqref{eq:WMAC}.
Given noise variance $\Exp{\abs{w_r}^2}$ and $(\dmax,\dmin)$, agent $r$ uses the $B'$ symbols to evaluate:
$y'_r := \frac{1}{B'}\sum_{b=1}^{B'} \abs{q'_r(b)}^2 - \Delta\ \Exp{\abs{w'_r}^2}$, where $\Delta := \dmax - \dmin$.
Then, for the $M$ sets of $B$ symbols, for each $m=1, \dots, M$, agent $r$ evaluates:
%\vspace{-0.75em}
\begin{align}\label{eq:postproc}
    y_r^{(m)} := \frac{\Delta}{B}\sum_{b=1}^B \abs{q_r^{(m)}(b)}^2 - \Delta\ \Exp{\abs{w_r}^2} + \dmin\ y'_r.
\end{align}
%\vspace{-1.25em}
%
It can be verified (see \cite[Lemma 1]{Agrawal2023}) that $y_r^{(m)}$ and $y'_r$ takes the following form:
\begin{align}\label{eq:otac_scalar}
    y_r^{(m)} = \sum_{j\in\sN_r} \nu_{jr}^{(m)} \lambda_j + \eta_r^{(m)}, \quad 
    y'_r = \sum_{j\in\sN_r} \nu'_{jr} + \eta'_r,
\end{align}
where, for all $j\in\sN_r$, we define $\nu_{jr}^{(m)} := \frac{P_j}{B} \sum_{b=1}^B \abs{\xi_{jr}^{(m)}(b)}^2$, and $\nu'_{jr} := \frac{P_j}{B'} \sum_{b'=1}^{B'} \abs{\xi'_{jr}(b')}^2$.
It turns out that the random variables $\eta_r^{(m)}$ and $\eta'_r$ are zero-mean, and with Assumption \ref{as:WMAC}(iii), we have also $\Exp{\nu_{jr}^{(m)}} = \Exp{\nu'_{jr}}$ for all $m$.

% Consequently, the scalar $y_r^{(m)}$ (similarly, $y'_r$, with $\lambda_k$ replaced with $1$) takes the form $y_r^{(m)} = \sum_{k\in\sN_r} c^{(m)}_{kr} \lambda_k + \eta^{(m)}_r$, where $c^{(m)}_{kr} := \frac{P_k}{B}\sum_{b=1}^B \abs{\xi^{(m)}_{kr}(b)}^2$ represents the random weight on the edge $(k,r)\in\sEi$ (see proof of Proposition \ref{prop:OTAC} for details).
% Note that $y'_r\to \sum_{k\in\sN_r} P_k \Exp{\abs{\xi_{kr}}^2}$ as $B'\to\infty$.
% Then, agent $r$ constructs the following two vectors:
% \begin{align}
%     \vy_r = \gamma_r\ [y_r^{(1)}, \dots, y_r^{(M)}]^T, \,\, \vy'_r = \gamma_r\ [y'_r, \dots, y'_r]^T \in\HS,
% \end{align}
% where $\gamma_r \in (0, (\Exp{y'_r})^{-1})$ is a design parameter.
% Using the information ($\vy_{r,i}$ and $\vy'_{r,i}$), 
%
The consensus protocol implemented by each $r\in\sN$ using information $(y_r^{(m)})$ and $y'_r$ obtained as above is: $(\forall  \iN)$
%\vspace{-0.5em}
\begin{align}\label{eq:agents_protocol}
     \gP_{r, i+1} = (\mI_M - \beta_i\ \diagm(\vy'_{r,i}))\ \gLki + \beta_i\ \vy_{r,i},
\end{align}
%\vspace{-1.5em}

\noindent where $\vy_r = \gamma_r\ [y_r^{(1)}, \dots, y_r^{(M)}]^T$, $\vy'_r = \gamma_r\ [y'_r, \dots, y'_r]^T \in\HS,$ and $\gamma_r \in (0, (\Exp{y'_r})^{-1})$ is a design parameter.
In practice, $\Exp{y'_r}$ can be estimated from past iterations.
% The protocol \eqref{eq:agents_protocol} can be rewritten in the form \eqref{eq:agent_scheme_consensus} by appropriately defining $\mPki$ and $\vnki$ as functions of random variables in the WMAC model \eqref{eq:WMAC} (see proof of Proposition \ref{prop:OTAC}).
%
The following proposition ensures that the sufficient conditions required by Theorem \ref{thm:T_QF} and Theorem \ref{thm:T_APSM} are satisfied by the proposed OTA-C protocol based consensus step in \eqref{eq:agents_protocol}.
%\vspace{-0.5em}
\begin{proposition}\label{prop:OTAC}
    Consider a system where agents exchange information using the proposed OTA-C protocol, where the conditions in Assumption \ref{as:WMAC} and Assumption \ref{as:net}(iii) are valid.
    Then, for every agent implementing the consensus step \eqref{eq:agents_protocol} in Scheme \eqref{eq:agent_scheme},
    the resulting communication model takes the form in \eqref{eq:comm_model},
    and the conditions in Assumption \ref{as:net}(i), \ref{as:net}(ii), and \ref{as:net}(iv) are satisfied.
    \begin{proof}
        Proof given in the Appendix \ref{prop:1}.
    \end{proof}
\end{proposition}
%\vspace{-0.5em}
\begin{remark}\label{rem:OTAC}
    Note that $B'$ can be chosen independently, which results in a reduced communication overhead compared to \cite{Agrawal2023}, where $B' = MB$.
    In addition, the requirement in \cite[Def.3.1]{Agrawal2023a} that every realization of graph $\sGi$ must produce a row-stochastic weight matrix is relaxed in this paper, leading to fewer constraints on the WMAC model.
    % Specifically, we relax the requirement in \cite[Def.3.1]{Agrawal2023a}, resulting in reduced communication overhead and fewer constraints on the WMAC.
    Moreover, we allow agents to independently select their transmit powers, which was previously unsupported in \cite{Agrawal2023,Agrawal2023a} as the graph was required to be undirected.
    %, as long as it can be ensured that the resulting graph remains strongly connected in expectation (Assumption \ref{as:net}(iii)).
    % Also, in contrast to \cite{Agrawal2023a}, support for directed graphs allows agents to independently select their transmit powers. %, as long as it can be ensured that the resulting graph remains strongly connected in expectation (Assumption \ref{as:net}(iii)).
    % Furthermore, when each agent $\kN$ selects the transmit power $P_k$ independently, the resulting graph is directed, which was previously not supported \cite{Agrawal2023a}.
    % Previously, in \cite{Agrawal2023}, the agents required to have equal transmit powers to obtain an undirected graph, necessitating coordination among the agents to agree on its value.
\end{remark}

%\vspace{-1em}
\section{Distributed machine learning application} \label{sec:sim}
%\vspace{-0.75em}
We simulate the task of supervised learning of a nonlinear function using data distributed over a decentralized network.
The data\footnote{The original dataset \cite{Seidov2020} is linearly scaled, and distances are in meters.}~consists of locations $\vx\in\sD:=[0, 1000]^3$ and corresponding measurements $f(\vx)\in[0,1]$.
At random times $(l_i)_\iN \subset\Natural$, agents move to a new location and obtain a noisy measurement $\ey_{k,l_i} = f(\vx_{k,l_i}) + e_{k,l_i}$, where $e_{k,l_i}$ is a zero-mean Gaussian random number with variance $0.09$.
At every $\iN$, each agent implements \eqref{eq:agent_scheme} using the sAPSM generator sequence $(\opTaki)_\iN$ in \eqref{eq:agent_scheme_local_opt} (as in Theorem \ref{thm:T_APSM}), and the OTA-C based consensus step \eqref{eq:agents_protocol}.
The sequence of bounded perturbations $(\zeta_i\vzki)_\iN$ is designed to promote sparsity in vector $\gLki$ which, as described later in this section, saves energy in communication using the proposed OTA-C protocol.
We only provide a brief description of the application here and refer the readers to \cite{Agrawal2023a} for more details.
% The explicit construction of mappings $(\opTaki)_\iN$, i.e., the cost functions $(\cT)_\iN$ and the sequence $(\zeta_i \vzki)_\iN$, is described in the following.

The cost functions $(\cT)_{\iN}$, for all $\kN$, are designed such that they satisfy the conditions in Assumption \ref{as:bound_sugGop}, and a solution of \eqref{eq:relax_set} gives a reasonable estimate of the function to be learned.
We use the multi-kernel approach \cite{Yukawa2012} with random Fourier features (RFF) approximations \cite{Rahimi2007} to model the nonlinear function $f$ as follows \cite{Shen2020}:
let $\hf(\vx) := \vh^T\ \pvv(\vx)$, where, for a design parameter $M>0,M\in\Natural$, the vector $\vh\in\Real^M$ is to be learned, and $\pvv:\sD\to\Real^M$ is the vector of RFF functions, fixed and known to all agents.
% $\hf(\vx) := \sum_{p\in\sP} \sum_{l\in\sL} h_{p,l}\ \pv_{p,l}(\vx) =: \vh^T\ \pvv(\vx),$
% where $\sP:=\{1,\dots,P\}$, $\sL:=\{1,\dots,L\}$, $\vx\in\sD$, 
% for each $l\in\sL$, the scalars $(h_{p,l})_{p\in\sP}$ are parameters to be learned, and the RFFs $(\pv_{p,l}(\vx))_{p\in\sP}$ (fixed and known to all agents) are generated apriori corresponding the kernel $\Ker_l$.
% Vectors $\vh$ and $\pvv(\vx)$ are the formed by stacking $PL$ scalars $(h_{p,l})$ and $(\pv_{p,l}(\vx))$, respectively, and we assume that $\vh\in\sX:= [-1, 1]^M$.
The cost function $\cT$ is defined as: 
$(\forall \iN)(\forall \kN)\ \cT(\vx) := \norm{\vx - \opP_{k,i}(\vx)} \norm{\gPki - \opP_{k,i}(\gPki)}$
where $\gPki$ is the estimate of agent $k$ at time $i$, and $\opP_{k,i}$ is the projection onto 
$\sQki := \{\vh : \abs{\vh^T \pvv(\vx_{k,l_i}) - \ey_{k,l_i}} \leq \varkappa_k\}$, where $l_i\in\acute{\Natural}$, $\varkappa_k\geq 0$ is a design parameter and 
$\vx_{k,l_i}\in\sD$ is the location of agent $k$ at time $i$.
In set $\sQki$, the parameter $\varkappa_k$ is chosen such that Assumption \ref{as:bound_sugGop} is satisfied with high probability.
The expression of projection mapping $\opP_{k,i}$ onto set $\sQki$ can be found in \cite{Yukawa2012}.
% Note that, in the proposed OTA-C protocol, the transmit signal of agent $k$ is zero whenever the scalar input $\lambda_k = \dmin$.
% We use this property to save energy in the prescribed system by promoting sparsity of $\gLki$ via reformulating the problem such that $\dmin=0$.
% For enforcing $\dmin=0$, we replace every $(h_{p,l})_{p,l}$ by two nonnegative parameters, one corresponding to $\pv_{p,l}$, and the other to $-\pv_{p,l}$.
% Hence, the total number of parameters in the reformulated setup is $2LP$. 

We assume that the set of optimal solutions is contained in $\sX := [\dmin,\dmax]^M$ where $\dmax=1$ and $\dmin=0$ for the sparsity-promoting scheme, otherwise $\dmin=-1$.
To get $\dmin=0$, model $\hf$ is modified by introducing another copy of the RFFs from the original model with a negative sign.
For the sparsity-promoting scheme, we design $\vzki$ as: $(\forall\kN)(\forall\iN)\ \vzki := \zeta_i^{-1}\ \left(\fq_{k,i}(\vy_{k,i}) - \vy_{k,i} \right)$,
where $\vy_{k,i} := \gPki - \opZki(\gPki)$ and $\fq_{k,i}:\Real^M\to\Real^M$ is defined element-wise, for each $m=1,\dots,M$, as
$\fq_{k,i}^{(m)}(x) := \text{sign}(x)\ \left[\abs{x} - \zeta_i(\abs{\vy_{k,i-1}[m]} + \varsigma_i)^{-1}\right]_+$,
where $\varsigma>0$ is a design parameter, $\text{sign}(a) := a/\abs{a}$, and $[a]_+ := \max(a, 0)$.
In essence, using the prescribed design, $(\zeta_i\vzki)_\iN$ reduces the reweighted $\ell_1$-norm \cite{Candes2008} of $\vy_{k,i} := \gPki - \opZki(\gPki)$, which leads to a sparse $\gLki$.
% With some simple modifications to the model $\hf$, agent $k$ do not need to transmit the zero components of $\gLki$, and hence, saving energy in communication.
% See \cite{Agrawal2023a} for details.
% where $\vzki$ is given by:

% \footnotetext{Sparsity inducing perturbations are designed to reduce the reweighted $\ell_1$-norm of the vector $\vy_{k,i}$, and operator $\fq$ is the proximal operator of the norm. See \cite{Candes2008,Agrawal2023a} for details.}

We simulate a time-varying network of $N=100$ agents as a geometric graph based on the locations sampled uniformly randomly from the dataset at random intervals.
The transmit power $P_k$ is sampled independently and randomly for each agent $\kN$ such that the resulting directed graph is strongly connected in expectation.
The channels $(\forall(k,r)\in\sE_i)\ \xi_{kr}$ and noise $(\forall r\in\sN)\ w_r$ are modeled as a circularly-symmetric zero-mean complex Gaussian random variables, with variance of the channel proportional to the inverse of squared distance (path loss), and variance of noise fixed to $-9dBm$ for all agents.
Agents randomly choose to either transmit or receive at any iteration $\iN$ (half-duplex system).
Other design variables are set to the following values: $M=50$, $B=20$, $B'=2B$, $(\forall \iN)\ \beta_i := (\lfloor i/50 \rfloor)^{-0.51}$, $\zeta_i := 10^{-6} (\lfloor i / 100 \rfloor + 1)^{-1}$, and $(\forall\iN)(\forall\kN)\ \rmuki = 0.5$.

We implement two schemes based on the proposed OTA-C protocols, and three standard communication protocols:
(OTAC-S): sparsity-promoting, over directed graphs;
(OTAC): without sparsity, over directed graphs;
% (OTAC-U): with sparsity, over undirected graphs;
(BDC): digital broadcast, where we assume Rayleigh fading channels with outage probability of $20\%$ at distance $500$m from transmitter;
% (APB): analog broadcast, where only one agent $k\in\sN$ broadcasts at any time $\iN$ over the fading channels $\xi_{kr}$, where we use $B=500$;
(NOC): no information sharing; and
(CEN): perfect (centralized and noiseless) information sharing.
Note that except OTAC, and BDC schemes, all other schemes introduce sparsity-promoting perturbations.

%\vspace{-0.5em}
\begin{figure}[!htb]
    \centering
    \begin{tikzpicture}
        \begin{axis}[
        scale only axis=true,
        width=0.8\columnwidth,
        height=0.6\columnwidth,
        label style={font=\small},
        tick label style={font=\small},
        xlabel={Iteration $i$},
        %axis y line*=left,
        ylabel={{\scriptsize NMSE (dB)}},
        mark size=1.5pt,
        mark repeat=10,
        grid=both,
        grid style=dashed,
        legend style={nodes={scale=0.8, transform shape}},
        ylabel near ticks,
        clip mode=individual,
        ]
        \addplot [solid, thick, mark=*] table[x=Iter, y=OTAC, col sep=comma]{plot_data_digraph_nmse.csv};
        \addplot [solid, thick, mark=square*] table[x=Iter, y=OTACS, col sep=comma]{plot_data_digraph_nmse.csv};
        \addplot [densely dashed, thick, mark=*] table[x=Iter, y=BDC, col sep=comma]{plot_data_digraph_nmse.csv};
        % \addplot [dashed, thick, mark=x] table[x=Iter, y=APB, col sep=comma]{plot_data_digraph_nmse.csv};
        \addplot [densely dotted, thick, mark=triangle*] table[x=Iter, y=NOC, col sep=comma]{plot_data_digraph_nmse.csv};
        \addplot [dotted, thick, mark=diamond*] table[x=Iter, y=CEN, col sep=comma]{plot_data_digraph_nmse.csv};
        \legend{OTACS, OTAC, BDC, NOC, CEN}
        \end{axis}
    \end{tikzpicture}
    \caption{Estimation error in terms of NMSE}
    \label{fig:nmse}
\end{figure}
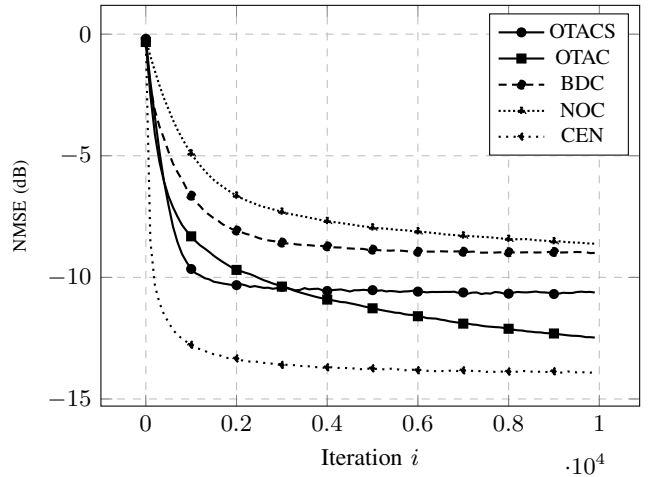
%\vspace{-0.5em}

The performance is compared in Figure \ref{fig:nmse} using a separate test dataset $\sR = \{(\vx,f(\vx))\}$ over the run of $10000$ iterations in time.
The metric used for comparison is the normalized mean square error (NMSE), given by 
\begin{align*}
    e(i) = \frac{1}{\abs{\sR}}\sum_{(\vx, f(\vx))\in\sR} \frac{1}{N} \sum_\kN \frac{\abs{\hf(\vhki;\vx) - f(\vx)}^2}{\abs{f(\vx)}^2},
\end{align*}
where each result is averaged over $100$ independent runs.
As expected, the best and worst performing scheme are CEN and NOC, respectively.
All three proposed OTA-C based schemes perform better than the BDC scheme, exhibiting the merits of OTA-C over the standard channel separation based strategy.
The OTAC scheme (without sparsity) gives the best results in the long run, whereas sparsity-promoting scheme OTAC-S shows faster convergence initially.
One plausible reasoning for this phenomenon is that proximity of the optimal estimate to the set of sparse vectors enables faster convergence in the beginning, but reaches an error floor as more information arrives.
% The OTAC-S scheme shows similar performance, but in directed-graph based OTAC-S, agents have the freedom to regulate their transmit powers independently.
It is worth noting that the sparsity-promoting schemes, i.e., OTAC-S, NOC and CEN, lead to sparse vectors $\gLki$ with less than $10\%$ nonzero entries, i.e., more than $80\%$ energy-saving in communication compared to OTAC and BDC schemes.

%In conclusion, this paper provides a unified framework for development and analysis of distributed algorithms that covers several (new and possibly yet-to-be invented) optimization and communication technologies.
%As a future extension, we are investigating methods for accelerating the convergence of our scheme by deriving theoretical bounds on the rate of convergence in one direction, and by developing algorithms that can capitalize on pioneering technologies such as beamforming in multi-antenna systems on the other.

\section{Conclusion}
The paper introduces a unified framework for the development and analysis of distributed algorithms, showcasing their adaptability to various optimization and communication technologies, both current and prospective.
Our future research will focus on deriving theoretical bounds on the convergence rate, with a goal towards faster solutions. 
Additionally, we aim to harness cutting-edge technologies like multi-antenna systems and low resolution ADC/DAC to further optimize and accelerate our algorithms in practice.
These exciting prospects promise to drive advancements in the field of distributed optimization and push the boundaries of what is achievable in dynamical systems.

\bibliographystyle{IEEEbib}
% \IEEEtriggeratref{13}
% \enlargethispage{-1in}
\bibliography{ICASSP_24_refs}

% \clearpage
\appendix
%!TeX root = proof_paper

We begin by reproducing a well-known result in the following, which we will use extensively in our proofs.
\begin{proposition}{(\cite[Theorem 1]{Robbins1971})}
    \label{prop:sto_app}
    Let $\pspace$ be the probability space, and $\sigmaF_1 \subset \sigmaF_2 \subset \dots$ be a sequence of sub-$\sigma$-algebras of $\sigmaF$.
    For each $n\in\Natural$, let $z_n$, $\beta_n$, $\xi_n$, and $\zeta_n$ be nonnegative $\sigmaF_n$-measurable random variables such that
    \begin{align}
        \label{eq:Robbins}
        \Exp{z_{n+1} \mid \sigmaF_n} \leq z_n (1 + \beta_n) - \zeta_n + \xi_n.
    \end{align}
    If, in addition, the series $\sum_{n\in\Natural}\beta_n$  and $\sum_{n\in\Natural} \xi_n$ converges $\Pas$, 
    then, the limit $\lim_{n\to\infty} z_n$ exists and finite, and the series $\sum_{n\in\Natural} \zeta_n$ converges.
\end{proposition}

\subsection{Proof of Theorem 1} \label{thm:1}
% \emph{Brief summary:}
% In a nutshell, the proof follows by designing a \emph{quadratic time-varying Lyapunov} (QTVL) function, as in \cite{Touri2013,Nedic2016}, and establishing conditions in Proposition \ref{prop:sto_app} for the sequence of scalars it generates.
% In this paper, in contrast to the results presented in \cite{Nedic2016}, new analysis and results are derived to account for the introduction of additive noise in the update dynamics \eqref{eq:agent_scheme}.
% % Particularly, the consensus step \eqref{eq:agent_scheme_consensus} presented in this paper involves a diminishing step-size sequence $(\beta_i)_\iN$, which necessitates new  to design the QTVL function leading to desired results (see Lemma \ref{lem:abs_seq} below).
% Although the QTVL function satisfies the properties of a norm at any time $i$, but, as the name suggests, it is a time-varying function.
% Hence, the result from \cite{Ermoliev1983} cannot be applied directly.
% Instead, the results of Theorem \ref{thm:T_QF} are derived directly from the conclusions of Proposition \ref{prop:sto_app}.

% \emph{Step-by-step proof:}
We begin by concatenating the consensus step \eqref{eq:agent_scheme_consensus} for all agents in the network $\sN$ as follows: $(\forall \iN)$
\begin{align}
    \gPii = ((1 - \beta_i) \mI_{MN} + \beta_i \mPi) \gLi + \beta_i\vni =: \mGb_i \gLi + \ve_i,
\end{align}
where, $\gPi$, $\vni$, and $\mPi$ are obtained by stacking $\gPki$, $\vnki$, and $\mPki$ column-wise for all agents $\kN$, respectively,
and recall that $\mI_{MN}$ is an identity matrix of size $MN$.
Here and henceforth, we define the following shorthand notations for convenience: $(\forall\iN)$
\begin{align*}
    &\mGb_i := \Exp{\mG_i}, \quad &&\mG_i := (1 - \beta_i) \mI_{MN} + \beta_i \mPi, \\
    &\ve_i := \beta_i(\mPt_i \gLi + \vni), \quad &&\mPt_i := \mP_i - \mPb_i, \quad \mPb_i := \Exp{\mPi} \\
    &\mGb_i = \mD_i \otimes \mI_M,  \quad &&\mD_i := (1-\beta_i)\mI_N + \beta_i \mA_i,
\end{align*}
where, in the last definition, we use the fact that $\mPb_i = \mA_i \otimes \mI_M$ (see discussion following equation \eqref{eq:comm_model}).
The matrix $\mGbki \in \Real^{M\times MN}$ is formed by taking $k$ rows of the matrix $\mGb_i$ corresponding to the $k$th agent, i.e., rows $(k-1)m+1$ to $km$ of the matrix $\mGb_i$.

% $\mGt_i := \mG_i - \mGb_i = \beta_i (\mPi - \mPb) =: \beta_i \mPt$, and
% Note that if $(\mA_i)_\iN$ satisfies the conditions of Assumption \ref{as:net}(i) and \ref{as:net}(iii), then the sequence of matrices $(\mD_i)_\iN$ also satisfies these conditions.
% Assumption \ref{as:net}(ii) is only satisfied for $\beta_i>0$.
% However, since $\lim_{i\to\infty} \beta_i = 0$, the condition in Assumption \ref{as:net}(ii) is not satisfied.

We will use the \emph{quadratic time-varying Lyapunov function} (described in the next paragraph) for our proofs. In the subsequent discussion, the following result from \cite{Touri2013,Nedic2016} will be helpful.
% The following lemma guarantees existence of a sequence of vectors $(\pi_i)_\iN$ in $\HS$ that are bounded away from zero, and satisfy the property \eqref{eq:abs_seq}.
\begin{lemma} \label{lem:abs_seq}
    Suppose that the sequence of matrices $(\mA_i)_\iN$ satisfies the conditions in Assumption \ref{as:net}(i)-(iii).
    Then, for the sequence of matrices $(\mD_i)_\iN$, where $(\forall \iN)\ \mD_i = (1-\beta_i)\mI_N + \beta_i\mA_i$, the following holds:
    \begin{align} \label{eq:abs_seq}
        (\forall \iN)\qquad \pi_i^T = \pi_{i+1}^T \mD_i,
    \end{align}
    where each vector $\pi_i\in\Realp^M$ is stochastic, i.e., sum up to one, and there exists some $\delta > 0$ such that $\pi_{k,i} \geq \delta$ for all $\kN$ and $\iN$, where $\pi_{k,i}$ is the $k$th coordinate of $\pi_i$.
    % Moreover, there exists some $\delta > 0$ such that $\pi_{k,i} \geq \delta$ for all $\kN$ and $\iN$, where $\pi_{k,i}$ is the $k$th coordinate of $\pi_i$.
    \begin{proof}
        The proof essentially follows from \cite[Lemma 9]{Touri2013} by establishing that the sequence $(\mD_i)$ satisfies the \emph{strong aperiodicity} and \emph{cut-balancedness} properties, and that each $\mD_i$ is a stochastic matrix.
        Then, by the definition of the class $\mathcal{P}^\star$ of sequence of matrices or \emph{chains} (cf.~\cite[Def.~3]{Touri2013}), the results of this Lemma follow immediately.
        
        Since $\mA_i$ is stochastic for all $\iN$ (Assumption \ref{as:net}(i)), the property that $\mD_i$ is stochastic follows from the definition of $\mD_i$.
        % Note that, since $\sum_{\iN} \beta_i$ diverges to infinity, for all $\iN$, there exists $\delta_0 > 0$ such that $\beta_i>\delta_0$.
        For a sequence of deterministic matrices, strong aperiodicity simply means that there exists $\gamma>0$ such that $[\mD_i]_{(p,p)} \geq \gamma $ for all $p\in\sN$ and $\iN$ \cite{Touri2013}.
        As $[\mA_i]_{(p,p)} \geq \epsilon > 0$, by definition, we have $[\mD_i]_{(p,p)} \geq (1 - \beta_i) + \beta_i \epsilon > \min(\epsilon, 1)$ since $\beta_i \in (0, 1)$.
        Hence, $(\mD_i)_\iN$ satisfies strong aperiodicity.

        Let $\sS$ be a nontrivial subset of agent indices, i.e., $\sS \subset \sN$ but $\sS \neq \sN$ or $\sS \neq \emptyset$, and define $\sSb := \sN \setminus \sS$ as complement of $\sS$.
        Then, since each $\mA_i$ is compliant with a strongly connected graph by Assumption \ref{as:net}(i) and \ref{as:net}(iii), there are edges with nonzero weights from a node in $\sS$ to $\sSb$ and vice versa.
        % Hence, for any two disjoint (nontrivial) set of indices $\sS$ and $\sSb$ in $\sN$, such that $\sS \cup \sSb = \sN$, there is an edge in either direction.
        It follows that, for all $\iN$, there exists a $\alpha > 0$ such that:
        \begin{align}\label{eq:cut-balance}
            \sum_{p\in\sS}\sum_{q\in\sSb} [\mA_i]_{(p,q)} \geq \alpha  \sum_{p\in\sS}\sum_{q\in\sSb} [\mA_i]_{(q,p)}.
        \end{align}
        Multiplying both sides of \eqref{eq:cut-balance} by $\beta_i > 0$, we get the same inequality for all $\mD_i$, since $[\mD_i]_{(p,q)} = \beta_i [\mA_i]_{(p,q)}$ for $p\neq q$.
        Hence, $(\mD_i)_\iN$ is cut-balanced with coefficient $\alpha > 0$.
    \end{proof}
\end{lemma}

Motivated by \cite{Touri2013,Nedic2016}, our proof uses a quadratic time-varying Lyapunov function $\opV(i, \vy)$, defined as: $(\forall \iN)(\vy\in\HS)$
\begin{align} \label{eq:funcV}
    \opV(i, \vy) := \sum_\kN \pi_{k, i} \norm{\gPki - \vy}^2,
\end{align}
where $\pi_i\in\Realpp^M$ is as specified in Lemma \ref{lem:abs_seq}.

In the following, we use the shorthand notation $\Expi{\cdot}$ to denote the conditional expectation $\Exp{\cdot \mid \gPi}$.
We proceed by verifying the conditions required for application of Proposition \ref{prop:sto_app} to the sequence $(\opV(i, \vy))_\iN$ for some $\vy\in\sQs$.
To this end, in the following, we first establish an inequality similar to \eqref{eq:Robbins} by bounding $\Expi{\opV(i+1, \vy)}$.
Incorporating the definition of $\opV$, we have:
$$\Expi{\opV(i+1, \vy)} = \sum_\kN \pikii \Expi{\norm{\gPkii - \vy}^2}.$$
% Using the definition of $\opTqki$ as in Definition \ref{def:T_QF}, we simplify the quantity $\Expi{\opV(i+1, \vy)}$ for some $\vy\in\sQs$ as follows:
% \begin{equation}\label{eq:orig_eq}
%     \begin{aligned}
%         &\Expi{\opV(i+1, \vy)} = \Expi{\sum_\kN \pikii \norm{\gPkii - \vy}^2} \\
%         &\quad = \sum_\kN \pikii \Expi{\norm{\opTqki(\gLki) - \vy}^2} \\
%         &\quad \leq \sum_\kN \pikii (\Expi{\norm{\gLki - \vy}^2} + \Expi{\epsilon_{k,i}})
%     \end{aligned}
% \end{equation}
% Inserting $\gLki$ from \eqref{eq:agent_scheme_consensus} into $\norm{\gLki - \vy}^2$ and expanding, we obtain:
The term $\Expi{\norm{\gPkii - \vy}^2}$ in the sum above can be further expanded by using \eqref{eq:agent_scheme_consensus} to replace $\gPkii$, as follows: $(\forall\kN)$
\begin{align}
    &\Expi{\norm{\gPkii - \vy}^2} = \Expi{\norm{\mGbki\gLi + \veki - \vy}^2} \nonumber\\
    &\quad = \Expi{\norm{\mGbki\gLi - \vy}^2} + \Expi{\norm{\veki}^2} \nonumber\\
    &\quad \quad + 2\Expi{\veki^T(\mGbki\gLi - \vy)} \nonumber\\
    &\quad = \Expi{\norm{\mGbki\gLi - \vy}^2} + \Expi{\norm{\veki}^2}, \label{eq:simp_1}
\end{align}
where the last expression follows from:
\begin{align*}
    &\Expi{\veki^T(\mGbki\gLi - \vy)} \\
    &\quad = \beta_i\Expi{(\mPtki \gLi + \vnki)}^T (\mGbki\gLi - \vy)\\
    &\quad = \beta_i(\Exp{\mPtki} \gLi + \Exp{\vnki})^T (\mGbki\gLi - \vy)\\
    &\quad = 0,
\end{align*}
where we use the properties that $\gLi$ is a (deterministic) function of $\gPi$, $\Exp{\mPtki}$ is zero by definition of $\mPtki$, and $\Exp{\vnki}$ is zero by Assumption \ref{as:net}(iv).

Before we proceed, we reproduce the following identity from \cite[Lemma 5]{Nedic2016}. 
For all $\vx, \vy$ in $\Real^N$, and scalar $c\in\Real$, where $\sum_n \vx_n = 1$, the following holds:
\begin{align} \label{eq:sq_id}
    (\vx^T\vy - c)^2 &= \sum_{n=1}^N \vx_n(\vy_n - c)^2 - \frac{1}{2} \sum_{n,n'=1}^N \vx_n \vx_{n'} (\vy_n - \vy_{n'})^2.
\end{align}

In the following, using the identity \eqref{eq:sq_id}, we would like to expand the expression $\sum_{\kN} \pikii \norm{\mGbki\gLi - \vy}^2$.
Expanding the expression of the norm $\norm{\mGbki\gLi - \vy}^2$ in \eqref{eq:simp_1}, defined for vectors in $\Real^M$, we obtain:
$$\norm{\mGbki\gLi - \vy}^2 = \sum_{m=1}^M ((\mGbkim)^T\gLi - \vym)^2,$$
where we define $\mGbkim\in\VS$ as the $f_{km} := (k-1)M + m$ row of $\mGb_i$.

In the following, we consider any $m$th summand $((\mGbkim)^T\gLi - \vym)^2$ in the above expression, and expand it as follows:
\begin{align*}
    &\sum_{\kN} \pikii ((\mGbkim)^T \gLi - \vym)^2 \\
    &\quad = \sum_{\kN} \pikii \Bigg[ \sum_{p\in\sN}\sum_{q\in\sM} \mGbkimp (\gLi^{(f_{pq})} - \vym)^2 \\
    &\quad - \frac{1}{2} \sum_{p,q} \sum_{p',q'} \mGbkimp \mGbkimpd (\gLi^{(f_{pq})} - \gLi^{(f_{p'q'})})^2 \Bigg] \\
    &\quad \mymark{(a)}{=} \sum_{\kN} \pikii \Bigg[\sum_{p\in\sN} [\mD_i]_{(k,p)} (\gL_{p,i}^{(m)} - \vym)^2 \\
    &\qquad - \frac{1}{2} \sum_{p\in\sN}\sum_{p'\in\sN} [\mD_i]_{(k,p)} [\mD_i]_{(k,p')} (\gL_{p,i}^{(m)} - \gL_{p',i}^{(m)})^2 \Bigg],
    % &\quad \mymark{(b)}{=} \sum_{p\in\sN} \pi_{p,i} (\gL_{p,i}^{(m)} - \vym)^2 \\
    % &\quad - \frac{1}{2} \sum_{\kN} \pikii \sum_{p\in\sN}\sum_{p'\in\sN} [\mD_i]_{(k,p)} [\mD_i]_{(k,p')} (\gL_{p,i}^{(m)} - \gL_{p',i}^{(m)})^2,
\end{align*}
where, for convenience, we define $\sM := \{1, \dots, M\}$. %, $f_{pq} := (p-1)M + q$, and use $\mGb_i^{(km,pq)}$ to denote the $(f_{km}, f_{pq})$ element of $\mGb_i$.
% Let us consider the element $\mGi(k,m,p,q) := [\mGi]_{((k-1)M+m, (p-1)M+q)}$ for a given tuple $((k,m),(p,q))$, where $k,p$ and $m,q$ take values in $\sN$ and $\sM:=\{1, \dots, M\}$, respectively.
Since matrix $\mGb_i = \mD \otimes \mI$, the element $\mGbkimp$ is nonzero only when $(p,k) \in \sEi$ and $m=q$, and in this case, $\mGbkimp = [\mD_i]_{(k,p)}$, i.e., the value is the same for all $m\in\sM$.
Therefore, by removing the zero elements in the summation $\sum_{p\in\sN}\sum_{q\in\sM} \mGbkimp$, we reduced it to $\sum_{p\in\sN} [\mD_i]_{(k,p)}$.
For the vector $\gLi$, the element $f_{(pq)}$ now corresponds to the element $f_{(pm)}$, since $q=m$, which in turn correspond to the $m$th element of $\gL_{p,i}$.
This leads to the expression $(a)$ above.
% The expression $(b)$ follows by recalling that $\mD_i$ satisfies \eqref{eq:abs_seq}, and changing the order of summations of the first term in the expression $(a)$ above.

In the first expression in the right-hand side of $(a)$ above, let $w_{p,i} := (\gL_{p,i} - y)^2$, where we suppressed the dependence on $m$ for convenience, and, with a slight abuse of notation, define $y := \vy^{(m)}$.
By defining $\vw_i$ as the vector stacking $w_{p,i}$ for all $p\in\sN$, we can write:
\begin{align*}
    \sum_{\kN} \pikii \left[\sum_{p\in\sN} [\mD_i]_{(k,p)} (\gL_{p,i}^{(m)} - \vym)^2\right] = \pi_{i+1}^T \mD_i \vw_i.
\end{align*}

% Using the definition of $\mD_i := (1 - \beta_i) \mI + \beta_i \mA_i$, and equation \eqref{eq:abs_seq}, we obtain:
% \begin{align*}
%     \pi_{i+1}^T \mD_i \vw_i = (1 - \beta_i) \pi_{i+1}^T \vw_i + \beta_i \pi_i^T \vw_i = \tpi_i^T\ \vw_i,
% \end{align*}
% where we define an auxiliary vector $\tpi_i := (1-\beta_i)\pi_{i+1}^T + \beta_i \pi_i^T$.
% Note that, since $(\forall\iN) \pi_i$ is a positive vector with elements bounded below by $\delta>0$ by Lemma \ref{lem:abs_seq}, we have that every element of $\tpi_i$, for all $\iN$, is bounded below by the same factor $\delta > 0$.

Using Lemma \ref{lem:abs_seq} above, we know that
\begin{align*}
    \pi_{i+1}^T \mD_i = \pi_i^T.
\end{align*}

In light of the discussion above, \eqref{eq:simp_1} leads to the following expression for $\Expi{\opV(i+1, \vy)}$:
\begin{align}\label{eq:temp1}
    \Expi{\opV(i+1, \vy)} = \sum_\kN \tpiki \norm{\gLki - \vy}^2 - a_i + \Tilde{b}_i,
\end{align}
where we define
\begin{align*}
    a_i &:= \frac{1}{2} \sum_{\kN} \pikii \sum_{p\in\sN}\sum_{p'\in\sN} [\mD_i]_{(k,p)} [\mD_i]_{(k,p')} \norm{\gL_{p,i} - \gL_{p',i}}^2, \\
    \Tilde{b}_i &:= \sum_\kN \pikii \Expi{\norm{\veki}^2}.
\end{align*}

Next, we replace $\gLki$ using \eqref{eq:agent_scheme_local_opt}, and using the definition of QFMS generator sequences in Definition \ref{def:T_QF}, it yields:
\begin{equation} \label{eq:opT_insert}
    \begin{aligned}
        \norm{\gLki - \vy}^2 &= \norm{\opTqki(\gPki) - \vy}^2 \\
        &\leq \norm{\gPki - \vy}^2 + \epsilon_{k,i}, \qquad \Pas
    \end{aligned}
\end{equation}

Therefore, by defining $b_i := \Tilde{b}_i + \sum_\kN \pikii \epsilon_{k,i}$, equation \eqref{eq:temp1} becomes:
\begin{align} \label{eq:thm1_cond}
    &\Expi{\opV(i+1, \vy)} \leq \opV(i, \vy) - a_i + b_i \qquad \Pas
\end{align}
Note that $a_i$ and $b_i$ are nonnegative random variables.
Moreover, $\sum_\iN b_i$ is convergent due to the following:

(a) $\Expi{\norm{\veki}^2} \leq \beta_i^2(\Exp{\norm{\mPtki^T\mPtki}_2}\Expi{\norm{\gLi}^2} + \Expi{\norm{\vnki}^2} + 2 \Expi{\norm{\mPtki^T \vnki} \norm{\gLi}})$, where $\Exp{\norm{\mPtki^T\mPtki}_2}$, $\Expi{\norm{\vnki}^2}$, and $\Expi{\norm{\mPtki^T \vnki}}$ are bounded ($\Pas$) by Assumption \ref{as:net}(iv), and $\Expi{\norm{\gLi}^2}$ is bounded ($\Pas$) because $\gLki$ is a vector in the compact set $\sX$ (see \eqref{eq:agent_scheme_local_opt}), and

(b) the series $\sum_\iN \epsilon_{k,i}$, for all $\kN$, and $\sum_\iN \beta_i^2$ are convergent by definition.

Applying Proposition \ref{prop:sto_app} on sequences $(\opV(i, \vy))_\iN$, $(a_i)_\iN$, and $(b_i)_\iN$, along with the inequality \eqref{eq:thm1_cond}, yields the following results:

(*) The sequence $(\opV(i, \vy))_\iN$ is convergent, $\Pas$, and

(**) the series $\sum_\iN a_i$ converges, $\Pas$

As shown later, convergence of the sequence $(\sum_\kN \tpiki \norm{\gPki - \vy}^2)_\iN$ established in (*) above, leads to the convergence of $(\norm{\gPki - \vy}^2)_\iN$.
However, we first need to establish result of part (ii) of Theorem \ref{thm:T_QF} for the proof of this statement.

% The result (*) gives that the sequence $(\sum_\kN \tpiki \norm{\gPki - \vy}^2)_\iN$ converges.
% Moreover, from the fact that $(\beta_i)_\iN$ converges to zero, we know that $\mD_i \to \mI$, and hence, $\pi_i \to \frac{1}{N}\vOne_N$.
% It follows that the sequence $(\sum_\kN \norm{\gPki - \vy}^2)_\iN = (\norm{\gPi - \vOne \otimes \vy}^2)_\iN$ converges $(\Pas)$.
% Furthermore, since $\tpiki$ is bounded below by $\delta > 0$, the sequence $(\norm{\gPki - \vy}^2)_\iN$ is bounded above.
% the series $\sum_\iN \sum_\kN \norm{\gPki - \vy}^2 = \sum_\iN \norm{\gPi - \vOne \otimes \vy}^2$ converges, $\Pas$
% Hence, for each $\kN$, the sequence $(\norm{\gPki - \vy}^2)_\iN$ converges, $\Pas$.
% And, since $\vy$ belongs to a bounded set $\sX$, the sequence $(\gPi)_\iN$ is bounded.
% This proves part (i) of Theorem \ref{thm:T_QF}.

From result (**), using the fact that $\piki$ is bounded below by $\delta > 0$ for all $\iN$ and $\kN$, we can lower bound $a_i$ as $a_i \leq (\delta/2) h_i$ for all $\iN$, where $h_i$ is given by:
% replacing $\mD_i$ in terms of $\mA_i$, we expand the expression of $a_i$ and use the property that $\pikii$ is bounded below, to obtain:
\begin{align*}
    &h_i := \sum_{(k,p,p')\in\sN^3} [\mD_i]_{(k,p)} [\mD_i]_{(k,p')} \norm{\gL_{p,i} - \gL_{p',i}}^2, \\
    &= 2 \sum_{p,p', p\neq p'} \beta_i (1-\beta_i) [\mA_i]_{(p,p')}  \norm{\gL_{p,i} - \gL_{p',i}}^2 \\
    &\quad + 2 \sum_{p,p', p\neq p'} \beta_i^2 [\mA_i]_{(p,p)} [\mA_i]_{(p,p')} \norm{\gL_{p,i} - \gL_{p',i}}^2 \\
    &\quad + \sum_{k, p, p', k\neq p \neq p'} \beta_i^2 [\mA_i]_{(k,p)} [\mA_i]_{(k,p')} \norm{\gL_{p,i} - \gL_{p',i}}^2.
\end{align*}
Thus, if the series $\sum_\iN a_i$ converges, then the series $\sum_\iN h_i$ must converge as well.
And for the series $\sum_\iN h_i$ to converge, each summand in the above expression of $h_i$ must converge.
It can be verified that the second and third summands converge due to the convergence of series $\sum_\iN \beta_i^2$, boundedness of $\gL_{p,i}$, and $[\mA_i]_{(p,p')}$ for all $(p,p')\in\sN^2$ and $\iN$.

Using that fact that for all $(p,q)\in\sE_i$, the corresponding component of $\mA_i$ is bounded below by $\epsilon$, i.e., $[\mA_i]_{(p,q)}\geq \epsilon$ (see Assumption \ref{as:net}),
the first summand can be lower bounded by $\epsilon \sum_{(p,p')\in\sE_i} \beta_i (1-\beta_i) \norm{\gL_{p,i} - \gL_{p',i}}^2$ ($\Pas$).
The Lemma \ref{lem:series} below shows that if $(\beta_i)_\iN$ is such that $\beta_i = (i+1)^{-\alpha}$ with $0.5 < \alpha \leq 1$, then
\begin{align}\label{eq:consensus_prop}
    \lim_{i\to\infty} \sum_{(p,p')\in\sE_i} \norm{\gL_{p,i} - \gL_{p',i}}^2 = 0, \quad \Pas
\end{align}
% This guarantees the existence of a subsequence $\lAN\subset \Natural$ such that for $l\in\lAN$ we have $\lim_{l\to\infty} \sum_{(p,p')\in\sE_l} \norm{\gL_{p,l} - \gL_{p',l}}^2 = 0$, $\Pas$
By Assumption \ref{as:net}(i) and (iii), i.e., strong connectivity of graphs $(\sG_l)_{l\in\lAN}$, we have that the limit must converge to zero for every pair of agents in the network.
More precisely, 
$$(\forall(p,q)\in\sN\times\sN)(\iN)\ \lim_{i\to\infty} \norm{\gL_{p,i} - \gL_{q,i}}^2 = 0, \ \ \Pas$$

\begin{lemma}\label{lem:series}
    Suppose that the series $\sum_\iN a_i b_i$ converges to some constant $c < \infty$, $\Pas$, where, for all $\iN$, $b_i \geq 0$, and $a_i = (i+1)^{-\alpha}$ with $0.5 < \alpha \leq 1$.
    Then, it follows that
    $$\lim_{i\to\infty} b_i = 0, \quad \Pas$$
    \begin{proof}
        We prove the statement by contradiction.
        There are two possible cases: (a) the sequence $(b_i)$ converges to some positive value $b > 0$, and (b) the sequence $(b_i)$ does not converge.
        We consider both cases and prove that each of them leads to a contradiction.

        \vspace{0.5em}

        \noindent \textit{Case (a)}: When $(b_i)_\iN$ converges to a constant $b > 0$.
        It follows that there exists some $T \in\Natural$ such that $\inf_{i>T} b_i = b$.
        In this case, $\sum_\iN a_i b_i > b \sum_\iN a_i = \infty$, which is a contradiction.

        \vspace{0.5em}

        \noindent \textit{Case (b)}: When $(b_i)_\iN$ does not converge.
        Note that by convergence of $\sum_\iN a_i b_i$ and the fact that $\sum_\iN a_i$ diverges to infinity, we know that $\liminf_{i\to\infty} b_i = 0$, i.e., there exists a subsequence $\lSN \subset \Natural$ for which $\lim_{l\to\infty, l\in\lSN} b_l = 0$.
        However, in this case, we can construct a subsequence $\Tilde{\Natural} \subset \Natural \setminus \lSN$, such that $\lim_{l\in \Tilde{\Natural}} b_l = \Tilde{b} > 0$.
        For any such subsequence $\Tilde{\Natural} \subset \Natural \setminus \lSN$, the series $\sum_{l\in\Tilde{\Natural}} a_l$ will still diverge to infinity.
        Hence, the series $\sum_{l\in\Tilde{\Natural}} a_l b_l > \Tilde{b} \sum_{l\in\Tilde{\Natural}} a_l = \infty$ diverges, which is again a contradiction.
    \end{proof}
\end{lemma}

In the following, we show that $\lim_{i\to\infty} \norm{\gL_{p,i} - \gL_{p',i}} = 0$ implies $\lim_{i\to\infty} \norm{\gP_{p,i} - \gP_{p',i}} = 0$.
Consider the sequence $(\abs{\gP_{p,i}^{(m)} - \gP_{p',i}^{(m)}})_\iN$, where $\gP^{(m)}$ denotes the $m$th coordinate of $\gP\in\HS$ for $m=1,\dots,M$.
The consensus step \eqref{eq:agent_scheme_consensus} corresponding to the $m$th coordinate is given by:
\begin{align*}
    &\abs{\gP_{p,i+1}^{(m)} - \gP_{p',i+1}^{(m)}} = \abs{(\mGb_{p, i}^{(m)} - \mGb_{p', i}^{(m)})^T \gLi + \ve_{p,i}^{(m)} - \ve_{p',i}^{(m)}} \\
    &\quad \leq \abs{(\mGb_{p, i}^{(m)} - \mGb_{p', i}^{(m)})^T \gLi} + \abs{\ve_{p,i}^{(m)} - \ve_{p',i}^{(m)}}.
\end{align*}

The following statements hold for each $m=1,\dots,M$, each $\iN$, and all $(p,p')\in\sN\times\sN$.
% Note that the sequence of time indices $i$ in the following discussion belong to the set $\lAN$, and hence, expressions like $\lim_{l\to\infty}$ must be read as the limit over the sequence of indices in the set $\lAN$.

(a) Since $(\beta_i)_\iN$ converges to zero, every coordinate of the sequence $(\ve_{k,i})_\iN$, where $\ve_{k,i} := \beta_i(\mPt_{k,i} \gL_i + \vn_{k,i})$, converges to zero $\Pas$
Consequently, for every $\delta_1>0$ there exists $N_1\in\iN$ such that for all $l \geq N_1$, we have $\abs{\ve_{p,l}^{(m)} - \ve_{p',l}^{(m)}} < \delta_1$ ($\Pas$).

% In the following, we show that for every $\delta_2>0$ there exists $N_2\in\Natural$ such that for all $i \geq N_2$, we have $\abs{(\mGb_{p, i}^{(m)} - \mGb_{p', i}^{(m)})^T \gLi} < \delta_2$.
% To this end, we use the property that $\lim_{i\to\infty} (a_i b_i)  = (\lim_{i\to\infty} a_i)(\lim_{i\to\infty} b_i)$.
(b) From the definition of $\mGb_i$, i.e., $\mGb_i := \mI_{MN} + \beta_i (\mPb_i - \mI)$, and the fact that $\lim_{i\to\infty} \beta_i = 0$, the limit $\lim_{i\to\infty} \mGb_i = \mI_{MN}$ exists.
This implies that the sequence $(\mGb_{p, i}^{(m)} - \mGb_{p', i}^{(m)})_\iN$ converges to the vector $\vu_{p,m} - \vu_{p', m}$, where $\vu_{p,m}\in\VS$ is a vector of all zeros except one at coordinate $(p-1)M + m$.
% Hence, $\lim_{i\to\infty} (\mGb_{p, i}^{(m)} - \mGb_{p', i}^{(m)}) \gLi = (\lim_{i\to\infty} (\mGb_{p, i}^{(m)} - \mGb_{p', i}^{(m)})^T)(\lim_{i\to\infty} \gLi) = $
Assume that $\gL^\star$ is a convergence point of the sequence $(\gL_{p,i})_\iN$.
Since $(\norm{\gL_{p,i}-\gL_{p',i}})_\iN$ converges to zero for all $(p,p')$, any convergence point $\gL^\star$ has the property that $(\vu_{p,m} - \vu_{p',m})^T \gL^\star = 0$, i.e., $(p-1)M + m$ and $(p'-1)M + m$ coordinates of $\gL^\star$ are the same for all $(p,p')$ and all $m$.
Using the property that, if $\lim_{i\to\infty} a_i$ and $\lim_{i\to\infty} b_i$ exists, then $\lim_{i\to\infty} (a_i b_i)  = (\lim_{i\to\infty} a_i)(\lim_{i\to\infty} b_i)$, we deduce that:
$\lim_{i\to\infty} (\mGb_{p, i}^{(m)} - \mGb_{p', i}^{(m)}) \gL_i = (\lim_{i\to\infty} (\mGb_{p, i}^{(m)} - \mGb_{p', i}^{(m)})^T) (\lim_{i\to\infty} \gL_i) = \vu_{p,m}^T\gL^\star - \vu_{p',m}^T\gL^\star = 0$.
Hence, for every $\delta_2>0$ there exists $N_2\in\iN$ such that for all $l \geq N_2$, we have $\abs{(\mGb_{p, l}^{(m)} - \mGb_{p', l}^{(m)})^T \gL_l} < \delta_2$.

Combining (a) and (b), we have the following result from \eqref{eq:agent_scheme_consensus}, that holds for each $m=1,\dots,M$ and all $(p,p')\in\sN^2$:
for every $\delta_3 >0$ there exists $N_3 = \max(N_1,N_2)$ such that for all $l\geq N_3$, we have $\abs{\gP_{p,l+1}^{(m)} - \gP_{p',l+1}^{(m)}} < \delta_3$.
It follows that $\lim_{l\to\infty} \norm{\gP_{p,l} - \gP_{p',l}} = 0$.
This proves part (ii) of Theorem \ref{thm:T_QF}.

In the following, we establish the convergence of $(\norm{\gP_{k,i} - \vy}^2)_\iN$ for all $\kN$.
% It can be verified that the sequence $(\tpi_i)_\iN$ converges to $\tps := (1/N) \vOne$.
% To see this, consider:
% \begin{align*}
%     \lim_{i\to\infty} \norm{\tpi_{i+1} - \tpi_i} = \lim_{i\to\infty} \beta_i \norm{\tpi_{i+1}^T (\mI_N - \mA_i)} = 0,
% \end{align*}
% due to the fact that $\lim_{i\to\infty} \beta_i = 0$.
% Moreover, 
We establish the convergence for each component $m=1,\dots,M$, and therefore, in the following, we only consider $m$th component of $\gP_{k,i} - \vy$, omitting the dependence on $m$ for convenience.
In other words, we only consider $\psi_{k,i} := \gP_{k,i}^{(m)}$ and $y := \vy^{(m)}$.
Let the limit in (*) be some constant $0 \leq c < \infty$, i.e., $\lim_{i\to\infty} \sum_\kN \tpi_{k,i} (\psi_{k,i} - y) = c$.
Define $v_{k,i} := \psi_{k,i} - y$ and $\vv_i = [v_{1,i},\dots,v_{N,i}]^T$ for the following discussion.
By definition of the limit, we have: for any $\epsilon_1>0$ there exists a $T_1\in\Natural$ such that for all $t>T_1$, it holds that $\abs{\tpi_t^T \vv_t - c}^2 < \epsilon_1$.
Also, from part (ii) of Theorem \ref{thm:T_QF}, it follows that: for any $\epsilon_2>0$ there exists a $T_2\in\lAN$ such that for all $t>T_2$, we have $\abs{\psi_{p,t} - \psi_{q,t}}^2 = \abs{v_{p,t} - v_{q,t}}^2 < \epsilon_2$ for all $(p,q)\in\sN\times\sN$.

Using \eqref{eq:sq_id}, we obtain the following for all $t\in\lAN$, $t > T_3 := \max(T_1, T_2)$:
\begin{align*}
    &\sum_\kN \tpi_{k,t} (v_{k,t} - c)^2 - \frac{1}{2}\sum_{p,q\in\sN^2} \tpi_{p,t} \tpi_{q,t} (v_{p,t} - v_{q,t})^2 < \epsilon_1, \\
    % &\Rightarrow \sum_\kN \tpi_{k,t} (v_{k,t} - c)^2 < \epsilon_1 + N^2 \epsilon_2, \\
    &\Rightarrow \sum_\kN (v_{k,t} - c)^2 < \frac{1}{\delta} (\epsilon_1 + N^2 \epsilon_2), \\
    &\Rightarrow \abs{v_{k,t} - c}^2 < \epsilon_3 := \frac{1}{\delta} (\epsilon_1 + N^2 \epsilon_2), \qquad \forall \kN,
\end{align*}
where the second expression follows from the fact that ${0 < \delta \leq \tpi_{p,t} \leq 1}$ for any $p\in\sN$ and $t\in\Natural$.
Since $\delta > 0$ is bounded away from zero, for any $\epsilon_3>0$, we can chose $T_3\in\lAN$ such that for all $t>T_3$ in $\lAN$, the inequality $\abs{v_{k,t} - c}^2 < \epsilon_3$ holds for every $\kN$.
This result is equivalent to $\lim_{t\to\infty} \abs{\psi_{k,t} - y}^2 = c < \infty$ for all $\kN$, i.e., every sequence $(\abs{\psi_{k,i} - y}^2)_\iNt$ converges $\Pas$
This proves part (i) of Theorem \ref{thm:T_QF}.

\begin{figure*}[htb]
    \begin{align*}
         \kappa_r(b) := \sum_{k\in\sN_r} \left[ \dmin P_k \left(\abs{\xi'_{kr}}^2 - \abs{\xi_{kr}}^2\right) + \left(\sum_{k\neq l} \xi_{kr}(b) \xi^*_{lr}(b) U_k(b) U^*_l(b) \sqrt{g_k(\lambda_k)g_l(\lambda_l)}\right) + w_r^*(b) \left(\xi_{kr}(b) s_k(b)\right)\right].
    \end{align*}
\end{figure*}

The proof of part \ref{thm:T_QF}(iii) follows the proof of \cite[Theorem 1(c)]{Ermoliev1983}.
By the assumption of Theorem \ref{thm:T_QF}(iii) that the sequence $(\sum_\kN \norm{\gPki - \gPs}^2)_\iNt$ converges, and from the results of Theorem \ref{thm:T_QF}(i), we have that the set of accumulation points of the sequence $(\gPki)_\iNt$ is nonempty ($\Pas$) for every $\kN$.
However, if $\sQs$ lies in a hyperplane, then there can be two accumulation points $\gPki'$ and $\gPki''$ equidistant from the hyperplane in which $\sQs$ lies.
In fact, any two points $\gPki'$ and $\gPki''$, with the property that $\norm{\gPki' - \gPs} = \norm{\gPki'' - \gPs}$ for any choice of $\gPs\in\sQs$, satisfy the result in \ref{thm:T_QF}(i).
However, by the additional condition specified in Theorem \ref{thm:T_QF}(iii), we ensure that $\sQs$ does not lie in a hyperplane, i.e., there exists a $u > 0$ such that, for any $\vs\in\sQs$, the set $\{\vh\in\HS \mid \norm{\vh - \vs} \leq u\}$ is nonempty.
Therefore, no two points $\gPki'$ and $\gPki''$  can exist that satisfy the property $\norm{\gPki' - \gPs} = \norm{\gPki'' - \gPs}$ for all $\gPs\in\sQs$.
Hence, the accumulation point must be unique for every $\kN$.
By the result \ref{thm:T_QF}(ii), we can also guarantee that this unique accumulation point is the same for all agents $\kN$, $\Pas$.
Hence, Theorem \ref{thm:T_QF}(iii) is proved.

\subsection{Proof of Theorem 2}  \label{thm:2}
Following the proof of Theorem \ref{thm:T_QF} above, we start analysis from \eqref{eq:opT_insert}:
\begin{equation}\label{eq:T_APSM_sim}
\begin{aligned}
    &\norm{\gLki - \vy}^2 = \norm{\opP_{\sX}(\gPki - \opZki(\gPki) + \zeta_i\vzki) - \vy}^2, \\
    &\quad\mymark{(a)}{\leq} \norm{\gPki - \opZki(\gPki) + \zeta_i\vzki - \vy}^2, \\
    &\quad = \norm{\gPki - \vy}^2 + \norm{\opZki(\gPki)}^2 - 2(\gPki - \vy)^T \opZki(\gPki) \\
    &\quad\qquad + \zeta_i^2 \norm{\vzki}^2 + 2\zeta_i \vzki^T(\gPki - \opZki(\gPki) - \vy).
\end{aligned}
\end{equation}
where $(a)$ follows from the property of orthogonal projection operators, namely, $(\forall \vx\in\HS, \forall \vy\in\sX)$ 
$\norm{\opP_{\sX}(\vx) - \vy}^2 \leq \norm{\vx - \vy}^2$.

Using the definition of $\opZki$ in \eqref{eq:subGop}, and method of proof in \cite[Theorem 2(a)]{Yamada2005},
the following relation immediately follows:
\begin{align}\label{eq:subG_sim}
    &\norm{\opZki(\gPki)}^2 - 2(\gPki - \vy)^T \opZki(\gPki) \notag\\
    &\qquad \leq - \rmuki (2 - \rmuki) \frac{(\cT(\gPki))^2}{\norm{\cS(\gPki)}^2}.
\end{align}

\sloppy Using \eqref{eq:subG_sim} in \eqref{eq:T_APSM_sim}, we obtain the same expression for $\Expi{\opV(i+1, \vy)}$ as in \eqref{eq:thm1_cond}, but with different definitions for $a_i$ and $b_i$, given by:
\begin{align*}
    a_i &:= \frac{1}{2} \sum_{\kN} \pikii \sum_{p\in\sN}\sum_{p'\in\sN} [\mD_i]_{(k,p)} [\mD_i]_{(k,p')} \norm{\gL_{p,i} - \gL_{p',i}}^2 \\
    &\quad + \sum_\kN \pikii \rmuki (2 - \rmuki) \frac{(\cT(\gPki))^2}{\norm{\cS(\gPki)}^2},
\end{align*}
\begin{align*}
    b_i &:= \sum_\kN \pikii \Expi{\norm{\veki}^2} + \zeta_i^2 \sum_\kN \pikii \Expi{\norm{\vzki}^2}  \\
    &\quad + 2\zeta_i \sum_\kN \pikii \Expi{\vzki^T(\gPki - \opZki(\gPki) - \vy)}.
\end{align*}

It can be verified that $a_i$ and $b_i$ are nonnegative for all $\iN$.
For proving that the series $\sum_\iN b_i$ converges, first note that the convergence of series $\sum_\kN \pikii \Expi{\norm{\veki}^2}$ is already established in the proof of Theorem \ref{thm:T_QF} above.
What remains is to prove that the second and third terms in the right-hand-side expression of $b_i$ are convergent.
In this direction, note that $\sum_\iN \zeta_i$ is convergent and there exists some $\kappa<\infty$ such that $\Expi{\norm{\vzki}}\leq\kappa$ for all $\iN$ and $\kN$ ($\Pas$) by definition (cf.~Definition \ref{def:T_APSM}).
By Cauchy-Schwartz inequality, we have $\vzki^T(\gPki - \opZki(\gPki) - \vy) \leq \norm{\vzki}\norm{\gPki - \opZki(\gPki) - \vy}$.
Hence, by Assumption \ref{as:bound_sugGop} and boundedness of set $\sX$ (vectors $\vy$ and $\gPki$ belong to $\sX$), there exists some $\varphi<\infty$ such that $\Expi{\norm{\vzki}\norm{\gPki - \opZki(\gPki) - \vy}} \leq \varphi$ for all $\iN$ and $\kN$, $\Pas$.

Therefore, as in the proof of Theorem \ref{thm:T_QF} above, applying the Proposition \ref{prop:sto_app}, in addition to the already established results in Theorem \ref{thm:T_QF}, we obtain the following result:
\begin{align*}
    (\forall \kN)\quad \lim_{i\to\infty} \frac{(\cT(\gPki))^2}{\norm{\cS(\gPki)}^2} = 0, \quad \Pas,
\end{align*}
where we used the convergence of series $(a_i)_\iN$, and the fact that $\pikii \rmuki (2 - \rmuki) > 0$ for all $\iN$ and $\kN$.
Since $\norm{\cS(\gPki)}$ is bounded by Assumption \ref{as:bound_sugGop}, we must have that $\lim_{i\to\infty} \cT(\gPki) = 0$ for all $\kN$, $\Pas$.
This proves part (i) of Theorem \ref{thm:T_APSM}.

Part (ii) of Theorem \ref{thm:T_APSM} can be proved by following the proof of \cite[Theorem 2(e)]{Cavalcante2011} line-by-line.

\subsection{Proof of Proposition 1}  \label{prop:1}
Note that the conditions in Assumption \ref{as:net}(i)-(iii) are satisfied by design.
An overview of the steps of proof is as follows:
First, we express the random variables $y_r^{(m)}$ and $y'_r$ in \eqref{eq:postproc} in terms of inputs $\lambda_k$ from agents in $\sN_r$.
Then, we transform \eqref{eq:agents_protocol} to the structure in \eqref{eq:agent_scheme_consensus} and obtain expressions for $\mPki$ and $\vnki$.
Using these expressions, we verify that $\mPki$ and $\vnki$ satisfy conditions in Assumption \ref{as:net}.

\begin{lemma} \label{lem:OTAC-rvs}
    Based on the OTC-protocol proposed in Section \ref{sec:comm}, the random variables $y_r^{(m)}$ and $y'_r$ in \eqref{eq:postproc} take the following form: $(\forall r\in\sN)(\forall m=1,\dots,M)$
    \begin{align}\label{eq:yy_exp}
        y_r^{(m)} := \sum_{k\in\sN_r} c^{(m)}_{kr} \lambda_k + \eta^{(m)}_r, \quad y'_r := \sum_{k\in\sN_r} c'_{kr} + \eta'_r,
    \end{align}
    where, we define
    \begin{align} \label{eq:edge_weights}
        c^{(m)}_{kr} := \frac{P_k}{B}\sum_{b=1}^B \abs{\xi^{(m)}_{kr}(b)}^2, \quad c'_{kr} := \frac{P_k}{B'} \sum_{b=1}^{B'} \abs{\xi'_{kr}(b)}^2,
    \end{align}
    and ${\eta^{(m)}_r := \frac{1}{B} \sum_{b=1}^B \big(\abs{w_r^{(m)}(b)}^2 - \Exp{\abs{w_r^{(m)}}^2} + \kappa_r(b)\big)}$, and the expression of $\kappa_r(b)$ is given at the top of this page, where we omitted $m$ for simplicity.
    The expression of $\eta'_r$ is obtained similarly by replacing $B$ with $B'$ and $\lambda_k$, for all $\kN$, in the expression of $\kappa_r(b)$.
    In addition, for all $k\in\sN_r$, the random variables $c^{(m)}_{kr}$ and $c'_{kr}$ have finite mean and variance,
    and the random variables $\eta^{(m)}_r$ and $\eta'_r$ are zero-mean and finite variance.
    Furthermore, there exists some $d<\infty$ such that $\Exp{\abs{\xi_{kr}}^2\ e_r} \leq d$ for all $r\in\sN$ and $k\in\sN_r$.
    % They represent the random weights over the edges of the graph $\sGi$, each corresponding to the $M+1$ realizations of the communication system.
    
    \begin{proof}
        Using definitions of the corresponding symbols, equation \eqref{eq:yy_exp} can be verified using elementary algebraic manipulations of the equation \eqref{eq:postproc}.
        Finite mean and variance of $c^{(m)}_{kr}$ and $c'_{kr}$ follows from Assumption \ref{as:WMAC}, where it is assumed that $\xi_{kr}$ has finite second and fourth moments.
        The noise terms $\eta^{(m)}_r$ and $\eta'_r$ are zero-mean because $e_r(b)$ has a mean zero.
        This follows from the fact that $U_k$ is zero-mean and i.i.d.
        The finite variance of $\eta^{(m)}_r$ and $\eta'_r$, again, follow from Assumption \ref{as:WMAC}, where, in addition to $\xi_{kr}$, it is assumed that $w_r$ also has finite second and fourth moments.
        Similarly, finiteness of $\Exp{\abs{\xi_{kr}}^2\ e_r}$ follows from finiteness of second and fourth moments of $\xi_{kr}$ and $w_r$.
    \end{proof}
\end{lemma}

Using expressions for $y_r$ and $y'_r$ from \eqref{eq:yy_exp}, the terms in agent $\kN$ consensus protocol \eqref{eq:agents_protocol} can be rearranged and stacked for all agents to obtain:
\begin{align*}
    \gP_{i+1} = (1 - \beta_i) \gLi + \beta_i (\mPi \gLi + \vni),
\end{align*}
where $\mPi$ is defined as follows: $(\forall (p,q)\in\sN\times\sN)(\forall m,n=1, \dots, M)$
\begin{align}\label{eq:new_comm}
    [\mP_i]_{(f_{pm},f_{qn})} = \begin{cases}
        \gamma_{p,i} c_{qp}^{(m)}, & \text{ if } (p,q)\in\sE_i, m = n, \\
        1 - \gamma_{p,i} \sum\limits_{k\in\sN_p} c'_{kp}, & \text{ if } p=q, m=n, \\
        0, & \text{otherwise}
    \end{cases}
\end{align}
where $f_{pm} = (p - 1)M + m$;
and elements of $\vni$ are given by: $(\forall p\in\sN)(\forall m=1,\dots,M)$
\begin{align} \label{eq:new_noise}
    \vni^{(f_{pm})} := \gamma_{p,i} (\eta^{(m)}_p - \eta'_r\gL_{p,i}^{(m)}).
\end{align}

It only remains to prove that $\mPi$ and $\vni$ defined above, satisfy the conditions in Assumption \ref{as:net}(iv).
To this end, noting that $\mPi$ and $\vni$ are functions of $(c^{(m)}_{kr}, c'_{kr})$ and $(\eta^{(m)}_r, \eta'_r)$, respectively, the conditions in Assumption \ref{as:net} follows from Lemma \ref{lem:OTAC-rvs}.

\end{document}